\documentclass[3p]{elsarticle}
\usepackage{hyperref}
\usepackage{amsmath,amssymb,amsthm}
\usepackage{latexsym}            % for the qed symbol

\usepackage{epsfig}
\usepackage{epstopdf}
\usepackage{amsfonts,amscd, makecell}
\usepackage[ruled,vlined]{algorithm2e}
\usepackage{color,comment}
\usepackage{amsbsy,bm}
\usepackage{subfigure}
\usepackage[utf8]{inputenc}          
\usepackage[T1]{fontenc}

\newcommand\vel{\mathbf{u}}

\newcommand\bv{\mathbf{v}}
\newcommand\bx{\boldsymbol{x}}

\newcommand\bn{\mathbf{n}}

\newcommand\bF{\mathbf{F}}

\newcommand\bD{\boldsymbol{D}}

\newcommand\btau{\boldsymbol{\tau}}
\newcommand\bsigma{\boldsymbol{\sigma}}

\newcommand{\be}{\begin{eqnarray}}
\newcommand{\ee}{\end{eqnarray}}
\newcommand{\ben}{\begin{eqnarray*}}
	\newcommand{\een}{\end{eqnarray*}}

\newcommand{\eps}{\varepsilon}
\newcommand {\na} {{\nabla}}
\newtheorem{theorem}{Theorem}[section]

\newtheorem{lemma}[theorem]{Lemma}

\allowdisplaybreaks \allowdisplaybreaks[4]

\hypersetup{
	colorlinks=true,       % false: boxed links; true: colored links
	linkcolor=black,          % color of internal links (change box color with linkbordercolor)
	citecolor=green,        % color of links to bibliography
	filecolor=magenta,      % color of file links
	urlcolor=cyan,         % color of external links
}

\begin{document}
	\begin{frontmatter}
		\title{Thermodynamically Consistent Diffuse Interface Model for Cell
			Adhesion and Aggregation}
		
		\author[mymainaddress]{Lingyue Shen}
		\author[mymainaddress]{Ping Lin}
	 \author[mysecondaryaddress]{Zhiliang Xu}	
		\author[mythirdaddress]{Shixin Xu  \corref{mycorrespondingauthor}}
		\cortext[mycorrespondingauthor]{Corresponding author}
		\ead{shixin.xu@dukekunshan.edu.cn}
		\address[mymainaddress]{Department of Mathematics, University of Dundee, Dundee DD1 4HN, Scotland, UnitedKingdom.}
		\address[mysecondaryaddress]{Department of Applied and Computational Mathematics and Statistics, University of Notre Dame,102G Crowley Hall, Notre Dame, IN 46556}
		\address[mythirdaddress]{Duke Kunshan University, 8 Duke Ave, Kunshan, Jiangsu, China. }
		\begin{abstract}	
			A thermodynamically consistent  phase-field model is introduced  for simulating multicellular deformation, aggreggation  under flow conditions.
			In particular,  a  Lennard-Jones type potential is proposed under phase field framework for cell-cell, cell-wall interactions. 
			A second-order accurate  in both space and time  $C^0$ finite element method is proposed to solved the model governing equations. Various
			numerical tests confirm the convergence, energy stability, and nonlinear mechanical properties of cells  of the proposed scheme.  Vesicles with different  adhesion  are also used to explain the  pathological risk  for patient with sickle cell disease. 
			
		\end{abstract}
		
	\end{frontmatter}

\section{Introduction}
Studying interaction between structures (cell-cell and cell-blood-vessel) is an important subject for understanding hemodynamics, because structural interactions at the cellular level unambiguously appear in a broad spectrum of blood flow related problems ranging from red blood cell (RBC) distribution \cite{Thomas2008aggregate} in blood vessel, the growth of blood clot \cite{fogelson2015fluid}, blood cell aggregation \cite{lee2016optical}, sickle cell disease \cite{alapan2014heterogeneous}, tumor cell dynamics \cite{connor2019mathematical} 
and diabetes  \cite{deng2020quantifying}.

Under the framework of sharp interface description in which interfaces separating different components of matter are idealized as hypersurfaces
with zero thickness, there are a number of studies focusing on  modeling  cell-cell and cell-blood-vessel interactions under blood flow conditions. Following the seminal work of Peskin \cite{PESKIN1972,peskin_2002},  the immersed boundary method was used to develop a model for platelet aggregation in \cite{fogelson2008immersed}. Cell-wall interaction model is introduced in \cite{connor2019mathematical} to simulate the adhesion and deformation of tumor cells at the vessel wall. Local and non-local models are described in \cite{gerisch2008mathematical} to investigate the invasion and growth of tumor cells. Cell-cell interaction modeling at the micro-scale is carried out in
\cite{ziherl2007flat,liu2006rheology,ziherl2007aggregates,zhang2008red}  to study the RBC aggregation problem. 
Multi-scale models are introduced in \cite{fedosov2014multiscale,xu2009study} in which dissipative particle dynamics (DPD) is used in \cite{fedosov2014multiscale} to establish a blood cell model in blood flow and a stochastic
cellular potts model (CPM) is introduced in \cite{xu2009study} for studying blood clot growth. In \cite{flormann2017buckling,quaife2019hydrodynamics,YNYang2019adhesion}, to account for cell-cell or cell-substrate adhesion, a Lennard-Jones type potential is introduced as a one dimensional function of the distance between the points on different cell membranes and substrates. The potential is a combination of a repulsive part and an attractive part which shows repulsion when particles get too close and shows attraction when the distance increase.
%\be
%\mu_{L-J}=\mathcal{H}\left[\left(\frac{\delta}{z}\right)^m-\frac a b \left(\frac{\delta}{z}\right)^n\right]~,
%\ee
%where $\mathcal{H}$ is the Hamaker constant, $z$ is the distance of the points, $\delta$ is the adhesion length scale and $m>n$ are integers that depend on the geometry and the molecular details of the two objects under adhesion. 
In \cite{quaife2019hydrodynamics}, this potential in combination with DPD method as mentioned before, is utilized to study cell deformation and doublet suspension.

In this work, we focus on utilizing diffuse interface method, also commonly known as phase-field approach \cite{AnMcFWh98} to model cell-cell and cell-blood-vessel interactions under blood flow conditions.  The phase-field models replace the sharp interface
description with a thin transition region in which microscopic mixing of the macroscopically distinct components of matter is permitted. The  phase-field approach not only yields systems of governing equations
that are better amenable to further analysis, but also handles topological changes of the interface naturally \cite{lowengrub_quasiincompressible_1998,AnMcFWh98,duqiang_2009_variational,duqiang_2005_curvature}. 
By treating the cell membrane as a diffuse interface, several phase-field vesicle interaction models  have been reported recently.  Fusion of cellular aggregates in biofabrication is modelled with a single phase-field function in \cite{yang2012modeling}. Vesicle-substrate adhesion model is established in  \cite{gu2014simulating,zhang2009phase,das2008adhesion}  where two phase-field functions are used to indicate wall and vesicle, respectively. Vesicle-vesicle adhesion model is introduced in \cite{gu2016two} in which different phase-field functions are used to represent different vesicles. Further, multiple cell aggregation model and  simulations are reported in \cite{jiang2019diffuse,wangxiaoqiang2016adhesion}.

%  \textcolor{blue}{Results of simulations using those model are in good agreement  with laboratory experiments  \cite{mills2004nonlinear,adhesion,Thomas2008aggregate}. }

Numerically simulating cell models with phase-field method is challenging. finite difference method is used in \cite{du2005phase} to study the elastic bending energy of the vesicle. Finite element method is introduced in  \cite{guillen2018unconditionally,zhang2009phase,duqiang_2008_bending,jiang2019diffuse}, spectral methods with an adequate number of Fourier modes is used in  \cite{gu2016two,cheng2018multiple}.However, few works developing energy stable numerical schemes are reported for vesicle or cell models, despite the fact that most phase-field models obey energy dissipation law at the continuous level. It is well-known that  preserving a discrete counterpart of the energy dissipation property of the system in numerical solution is crucial for long term stability of numerical computation   \cite{hua2011energy,lin2006simulations}. In this paper, we also try to develop a numerical scheme retaining discrete energy dissipation law that is consistent with the continuous one. 

%Layer-by-layer decomposition is introduced in \cite{yang2012modeling}, finite element method is used in \cite{jiang2019diffuse}. A Spectral method with an adequate number of Fourier modes is used in \cite{gu2016two}.

%In numerical simulation, the study are similar as the ones mentioned in Section \ref{review2}. Lattice Bolzmann method and immersed finite element method are reported in \cite{zhang2008red} and \cite{liu2006rheology} respectively under sharp interface condition. For diffuse interface model, layer-by-layer decomposition is used in \cite{yang2012modeling}. However, few of the works are reported along with an energy stable discrete scheme despite energy decaying law in continuous form. 
%In the paper, an energy stable $C^0$ finite element scheme is established based on a Chan-Hilliard type interface model of multiple vesicle. Energy variational method is used for continuous model derivation. simulation results indicates that the model fits the experiment well.

The purposes of this paper is twofold.
The first goal of this paper is to derive a thermodynamically consistent phase-field model for vesicles' motion, shape transformation   and their structural interactions under flow conditions in a closed spatial domain using an energy variational method \cite{shen2020energy,guo2015thermodynamically}.  All the physics that ones are interested in are taken into consideration  through  definitions of the energy functional and the dissipation functional, together with the kinematic relations of dynamic  evolution of model state variables. 

In this work, all terms of the energy and dissipation functionals are defined on the bulk region of the domain  including cell-cell adhesion and cell-wall interactions.   A  Lennard-Jones typ potential within the framework of phase-field method is derived for the vesicle-vesicle interaction, in which an interacting potential including both repulsive and attractive parts is established with respect to the phase-field order parameter. To consider the cell-wall interaction, the wall is just treated as 
a new phase-field order parameter or function $\phi_w(\bx)$   in the proposed Lennard-Jones type potential. Then, the cell phase-field function $\phi$ defined on the bulk region of the domain which would be used in the derivation of model governing equations by the energy variational approach \cite{shen2022energy}. The advantage of doing so is that the variational procedure only needs to consider the bulk unlike previous works \cite{shen2020energy,shen2022energy} where the boundary needs to be considered since the wall effect are included in boundary conditions of the system, which simplifies the model derivation considerably.   %\textcolor{red}{Move to the 3rd paragraph???} A few works have been published about modeling of the interaction between cells\cite{wangxiaoqiang2016adhesion,YNYang2019adhesion}. Lennard-Jones potential is applied to model the interaction in \cite{YNYang2019adhesion}, an energy related to the inter-cell attraction is introduced using phase-field method in  With phase-field theory in \cite{wangxiaoqiang2016adhesion}. Results of simulations using those model are in good agreement  with laboratory experiments  \cite{mills2004nonlinear,adhesion,Thomas2008aggregate}. 

%To account for the  boundary effect in the new model,  a new phase parameter for the domain boundary is defined. \textcolor{red}{please elaborate}   

Performing variation of the energy and dissipation functionals yields a Cahn-Hilliard-Navier-Stokes (CH-NS) system with Allen-Cahn (AC) general Navier boundary conditions (GNBC) \cite{qian_2006_slipBC,shen2022energy}. 
This is in contrast to  most previous works \cite{duqiang_2004_bending, duqiang_2008_bending, chen2015decoupled} in which boundary effects were  rarely derived during the course of model derivation. Dirichlet or Neumann type boundary conditions were simply added to these models at the end to close the governing equations \cite{voigt_2014_local_inextensibility,duqiang_2009_variational,duqiang_2005_curvature,hua2011energy,guo2021diffuse}.

Moreover,  our model  accounts for the incompressibility of the fluid, and the local and global inextensibility of the vesicle membrane  by introducing two Lagrangian multipliers, hydrostatic pressure $p$ and surface pressure $\lambda$ \cite{laimingzhi_2020_immersed_boundary} and penalty terms, respectively. 

The second goal of this paper is to design an energy stable finite element scheme for solving model governing equations. This scheme preserves a discrete counterpart of the energy dissipation property of the model in numerical solution, which is crucial for long term stability of numerical computation. 
We note that although this scheme is designed for solving the model equations, it can be easily applied to other CH-NS systems.    We apply our model to simulate multiple vesicles interacting with each other and the domain boundaries under flow conditions \textit{in silico}, such as red blood cells passing a branched blood vessel. \textcolor{red}{elaborate}

Rest of the paper is organized as follows.  In Section \ref{sec:model}, the thermodynamically consistent model considering cell-cell and cell-wall interaction is derived,   the  energy dissipation law of the model is given. Then the discrete scheme of the model and the discrete energy law that is consistent with the continuous condition is given in Section \ref{sec:numericalscheme}. Section \ref{sec:results} is used to present results of numerical simulations  and compare them with the data collected in laboratory experiments.  The conclusion is drawn in Section \ref{sec:conclusion}.

\section{ Model  Derivation} 
\label{sec:model}
Derivation of the model in this paper is based on the energy dissipation law which holds ubiquitously in  physical systems involving irreversible processes \cite{xu2019three,eisenberg2010energy,hyon2010energetic,xu2014energetic,guo2021diffuse}. This law states that for an isothermal and closed system the  rate of change of the energy of the system is equal to the dissipation of the energy as  follows: 
\begin{equation}
\label{ED_1}
\frac{d}{d t} E^{\rm total} = - \triangle \leq 0~,
\end{equation}
where $E^{\rm total}$ is the total energy of the system, which is the sum of the kinetic energy $\mathcal{K}$ and the Helmholtz free energy $\mathcal{F}$ of the system, and $\triangle$ is the rate of energy dissipation, which in fact is entropy production. Eq.~(\ref{ED_1}) can be easily derived via the combination of
the First and Second Laws of Thermodynamics. The choices of the total energy functional and the dissipation functional, together with the kinematic
(transport) relations of the variables employed in the system, determine all the physical
and mechanical considerations and assumptions for the problem \cite{Giga2017VariationalMA}. A general technique to determine these relations of variables may be found in \cite{guo2015thermodynamically,guo2021diffuse}. 

% \textcolor{blue}{We assume that $\triangle$ is induced by viscosity, diffusion, membrane surface relaxation and friction between cell and the wall. move the other place???} 

\subsection{  Multi-vesicle interaction System}
Energy variational method \cite{shen2020energy} is adopted for model derivation to ensure the thermodynamics consistency of the derived model. 
Specifically, the model respects the energy dissipation law given by (\ref{ED_1}). The model derivation begins with defining functionals of the total energy and dissipation of the system, respectively,  and making the kinematic assumptions of the state variables of the system based on physics laws of conservation. We refer the readers to  \cite{shen2020energy,shen2022energy,xu2018osmosis} for detailed discussions of this method.

\begin{figure}[!ht]
	\centering
	\includegraphics[width=3.5in]{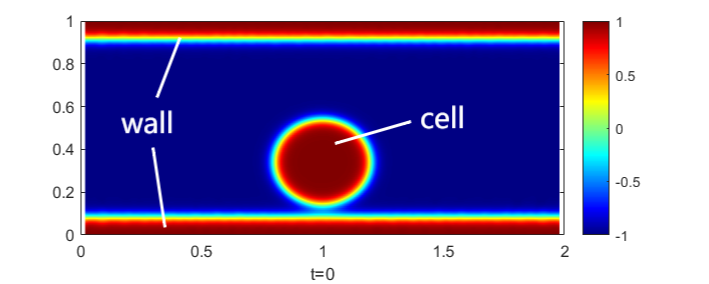}
	\caption{ A schematic showing the model setup. Phase-field functions are used to represent the domain wall boundary and the cell interface. }
	\label{fig:set up}
\end{figure}

Figure \ref{fig:set up} is a schematic showing the setup of the model. Let the problem domain be $\Omega$, and its wall boundary be $\partial \Omega_w$. For a multi-cell flow system, the dynamic evolution of the $i^{th}$ cell (or vesicle) under flow conditions within $\Omega$ is tracked by the phase-field function  $\phi_i(\boldsymbol x, t)$. Notice that $\phi_i(\boldsymbol x, t) \in [-1, 1]$ with  $\phi_i=1$ for intracellular space and $\phi_i=-1$ for extracellular space. The membrane of the cell is identified by $\phi_i(\boldsymbol x, t)=0$.  We also introduce a  phase-field function  $\phi_w(\boldsymbol x) \in [-1,1]$ to represent the wall boundary of the domain as shown in Fig.~\ref{fig:set up}. This is for considering cell-wall interaction described below. 
%We introduce a new phase $\phi_w$ to describe  the wall boundary  as shown in Figure \ref{fig:init}. The  thickness of the phase $\phi_w$ is chosen to be $\Gamma_{wall}=2\times 10^{-7}m$ \cite{gu2014simulating}. This is used for modeling cell-wall interaction...

We assume that the dynamical evolution of the phase-field function $\phi_i$ is a generalized gradient flow,  see Eq.~\eqref{assumption1}.
We also utilize the  
laws of conservation for describing the dynamics of the  momentum and the total mass of the system, see Eqs.~ \eqref{assumption2}-\eqref{assumption3}. The cell membrane (interface) is assumed to be inextensible. The equation for  local inextensibility of the interface  is given by Eq.~\eqref{assumption4}.
\begin{subequations}
	\label{assumption}
	\begin{eqnarray} 
	\frac{\partial \phi_i}{\partial t} +\nabla\cdot (\vel \phi_i) = \nabla\cdot {\mathbf q}_{\phi_i}~, \label{assumption1}\\
	\rho(\frac{\partial\vel}{\partial t}+(\vel\cdot\nabla)\vel) = \nabla\cdot\bsigma_{\eta}+\bF_{\phi_1, \phi_2,..., \phi_N}~,\label{assumption2}\\
	\nabla \cdot\vel = 0~,\label{assumption3}\\
	(\mathcal{P}_i: \na\vel)\delta_i=0 ~,\label{assumption4} 
	\end{eqnarray}
\end{subequations}
where  $\rho$ is the averaged density of the system. In this work,   $\rho$  is assumed to be a constant. $\vel$  is the macroscopic velocity of the system,  and  $\bsigma_{\eta}$ is the system's visco-elastic stress to be determined.

The  undetermined flux $\mathbf{q}_{\phi_i}$ in Eq.~(\ref{assumption1}) will be determined by postulating that  $\phi_i$ is driven by gradients in a chemical potential. This leads to the Cahn-Hilliard equation which ensures  conservation of the  volume of a cell during its dynamic evolution. $\boldsymbol{F}_{\phi_1,\phi_2,...,\phi_n}$ in Eq.~(\ref{assumption2}) is the  body force induced by vesicle-fluid interaction and is to be determined as well. 

Equation (\ref{assumption4}) is the diffuse interface approximation of the  local inextensibility  of the membrane of the $i^{th}$ cell  with relaxation \cite{shen2022energy,voigt_2014_local_inextensibility}.   
$\delta_i = \frac 1 2 \gamma^2|\nabla\phi_i|^2$ is the surface delta function with the diffuse interface thickness $\gamma$. $\mathcal{P}_i$ is the projection operator, and is  defined to be $(I-\bn_{i,m}\otimes \bn_{i,m})$. $\bn_{i, m}=\frac{\nabla\phi_i}{|\nabla\phi_i|}$ is the unit outward normal vector of the interface when it is defined as an implicit surface by the phase-field function.  This equation is equivalent to $\mathcal{P}_i:\na\vel=0$ under a sharp interface condition. In phase field frame, this local inextensibility is extended to the whole domain by multiplying a scalar function $\delta_i$ \cite{shen2020energy} for the convenience of computation. Following the idea in \cite{voigt_2014_local_inextensibility,shen2022energy}, here we add a relaxation term $\xi\gamma^2\na\dot(\phi_i^2\na\lambda_i)$ in Eq. \eqref{assumption4}  where $\xi$ is a parameter independent of $\gamma$,  $\lambda_i$ is a function that measures the interface “pressure” induced by the inextensibility of the membrane of the $i^{th}$ cell.  Thus Eq. \eqref{assumption4} now takes the form:
\be\label{local_inex}
(\mathcal{P}_i: \na\vel)\delta_i+\xi\gamma^2\na\dot(\phi_i^2\na\lambda_i)=0
\ee
Eq.\eqref{assumption1}-\eqref{assumption3} and \eqref{local_inex} together constitute the governing equations of the system. 

The boundary conditions on the top and bottom of the domain, denoted by $\partial\Omega_w$, are given as follows:
\be
\label{assumption_bd}
\left\{\begin{array}{l}
	\vel\cdot\bn = 0~,  \\ 
	\vel_{\tau}\cdot \btau_k = f_{\btau_k}~, \\
	\frac{\partial \phi_i}{\partial t} + \vel\cdot\nabla_{\Gamma}\phi_i = J_{\Gamma_i}~, \\
	\mathbf{q}_{\phi_i} \cdot \bn = 0~ .
	%		\frac{\partial f}{\partial t} + \vel\cdot\nabla_{\Gamma}f = -J_{f}, 
\end{array}\right.
\ee
%where $(\cdot)_i, i=1,2,3,...N$ refers to the \textcolor{blue}{corresponded physical quantities of different cells like phase, chemical potential...} \textcolor{red}{what is $(\cdot)_i$?}.

On the boundary $\partial\Omega_w$, an Allen-Cahn $(\ref{assumption_bd})_3$  type boundary condition is employed for $\phi_i$. $\nabla_{\Gamma} = \nabla-\bn(\bn\cdot\nabla)$ is the surface gradient operator, and $\vel_{\tau} = \vel-(\vel\cdot\bn)\bn$ is the fluid slip velocity with respect to the wall where $\btau_i, i=1,2$ are the tangential directions of the wall surface. 
$\bn$ is the unit outward normal vector of the wall.    $f_{\btau_i}$ is the slip velocity of the fluid on the wall along the $\btau_i$ direction.  $J_\Gamma$ represents the Allan-Cahn type of relaxation on the wall by using the phase-field method, and is to be determined, see \cite{qian_variational_2006,shen2020energy}.   

The total energy functional $E_{total}$ for the multi-cell system  is assumed to be the sum of the kinetic energy $E_{kin}$ in the macroscale, elastic energy $E_{cell}$ of cell membrane, cell-cell interaction energy  $E_{int}$ and cell-wall adhesion energy   $E_{w}$ in the microscale. Therefore  
\be\label{energy}
E_{total}&=& \underbrace{E_{kin}}_{Macroscale}+ \underbrace{E_{cell}+E_{int}+E_{w}}_{Microscale}~\\
&=&\int_{\Omega}\left(\frac 1 2 \rho | \vel|^2\right)d \boldsymbol{x}+\sum_{i=1}^{N}\int_{\Omega} \frac{\hat\kappa_B}{2\gamma}\left|\frac{f(\phi_i)}{\gamma}\right|^2d\bx +\sum_{i=1}^{N}\frac {\mathcal{M}_s} {2} \frac {\left(S(\phi_i(t))-S(\phi_{i}(t=0))\right)^2} {S(\phi_{i}(t=0))}\nonumber\\
&+& \int_{\Omega} Hd \boldsymbol{x}+\sum_{i=1}^{N}\int_{\Omega}f_w(\phi_i)ds~.
\ee
%Here $\rho$ is the density of the fluid, and  is taken to be a constant in this work. $\vel$ is the velocity.
Here $\hat\kappa_B$ is the bending modulus of the  membrane, and $\gamma$ is the thickness of the diffuse interface. The membrane energy density is given by
\be
\label{free energy phi}
f(\phi_i)&=&\frac{\delta G}{\delta\phi_i} = -\gamma^2 \Delta\phi_i+ (\phi_i^2-1)\phi_i~,
\ee
with 
\be
G(\phi_i)&=& \frac{\gamma^2|\nabla\phi_i|^2}2+\frac{(1-\phi_i^2)^2}{4}~.
\ee
The function $S(\phi_i) = \int_{\Omega}\frac{G(\phi_i)}{\gamma}dx$ is used to measure the surface area of the cell \cite{voigt_2014_local_inextensibility,duqiang_2009_variational,shen2020energy}. $\mathcal{M}_s$ is the penalty constant. 

The term $H$ denotes the  interaction energy density induced by the interaction of cells.    There are many different ways to define $H$ \cite{gu2016two,jiang2019diffuse,yang2012modeling}.  Here we begin with considering  mechanical interaction between two cells identified by phase-field functions $\phi_1$ and $\phi_2$, respectively.  Recall that $\phi_i=1$ represents the intracellular space, and $\phi_i=-1$ represents the extracellular space of the $i^{\rm th}$ cell, respectively. 
Whether there exists  mechanical interaction between the two cells can be determined by measuring the overlapping (i.e., occupying the same physical space) of the membranes of these two cells.  The following two functions are introduced to measure the overlapping  
$$d_1(\bx)=(\phi_1(\boldsymbol x)-1)(\phi_2(\boldsymbol x)-1)~,~~ \boldsymbol x\in \Omega~,$$ and  
$$d_2(\bx)=(\phi_1(\boldsymbol x)+1)(\phi_2(\boldsymbol x)+1)~,~~ \boldsymbol x\in \Omega~.$$ 

It is obvious that $d_1, d_2\in[0,4]$. From physical point of view, we can regard $d_1$ as the ``distance'' in the space of phase-field functions since it reaches maximum when both $\phi_1$ and $\phi_2$ equal to -1 which means there is no overlap between the two cells. $d_1$ reaches minimum when two cells fully overlap, i.e., $\phi_1=\phi_2=1$. Similarly, $d_2$ could be regarded as the extent of overlapping since it reaches maximum when cells fully overlap, and reaches minimum when there is no overlap.  Then we propose the interaction potential $H$ by a polynomial function of $d_2$ and $d_1$:
\be
H=Q(d_2^{\alpha}+C)d_1^{\beta}=Q C(\phi_1+1)^{\beta}(\phi_2+1)^{\beta}-Q(\phi_1^2-1)^{\alpha}(\phi_2-1)^{\alpha}(\phi_1+1)^{\beta-\alpha}(\phi_2+1)^{\beta-\alpha}~,
\label{eq:interact_engy}
\ee
where $Q$ and $C$ are constant that controls the minimum of the interaction energy. 
The first term of Eq.~(\ref{eq:interact_engy}) accounts for the repulsion between the two cells  which increase sharply to prevent fully overlapping of the cells, and the second terms is for the attraction between cells  which only appears when the cells contact each other, which means the they start to overlap.  

In this work, we take $\alpha=\beta=2$. Thus  the interaction potential (\ref{eq:interact_engy}) becomes:
\be\label{H2phase}
H= Q_1  (\phi_1+1)^2(\phi_2+1)^2- Q_2 (\phi_1^2-1)^2(\phi_2^2-1)^2~.
\ee

\begin{figure}[!ht]
	\centering
	\includegraphics[width=3.5in]{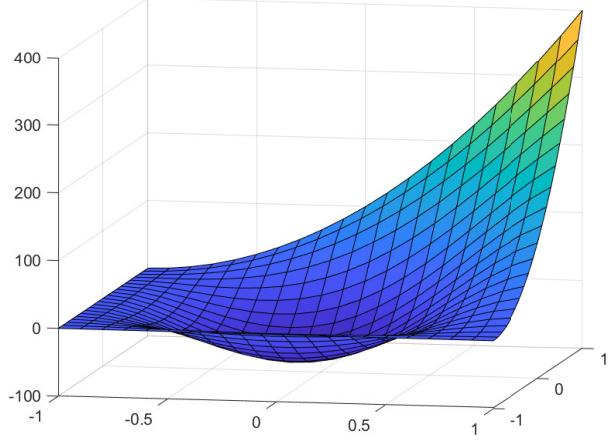}
	\caption{Interaction energy with respect to $\phi_1$ and $\phi_2$ with $Q_1=50$, $Q_2=400$. }
	\label{fig:interaction energy}
\end{figure}

Now we may consider $H$ as a 2D function of $\phi_1$ and $\phi_2$, where both $\phi_1$ and $\phi_2$ take values between -1 and 1 and their values represent relative position of vesicle 1 and vesicle 2 when they overlap. Also we need to point out that it is because $H$ is regarded as a 2D function, we need two variables ($d_1, d_2$) or ($\phi_1$, $\phi_2$) to formulate a single-concave potential energy with only one minima in the domain. Figure \ref{fig:interaction energy} plots the energy landscape of the interaction potential due to presence of the phases $\phi_1$ and $\phi_2$.   The energy is equal to 0 when two phases  do  not touch or overlap i.e., $\phi_1=\phi_2=-1$. When they start to overlap, the energy firstly decreases which means that the attraction force between these two phases dominates. Then the energy goes up which indicates the repulsive force dominates. This prevents the two phases from occupying the same physical space.  So the interacting potential energy behaves conceptually similar to a 2D kind of Lennard-Jones potential. We should point out that Eq. \eqref{eq:interact_engy} or \eqref{H2phase} is just constructed to mimic Lennard-Jones type of repulsive-attractive feature. Other formula with similar behavior may be good too. We remark that in \cite{gu2014simulating} a potential of form $(\phi_1^2-1)^2 (\phi_2^2 -1)$ is used, where only attractive feature is included.

More general, the interaction energy in multiple phases case goes:
\be\label{potential energy}
H=\sum_{i=1}^N\sum_{i<j}H_{ij}
=\sum_{i=1}^N\sum_{i<j} \left[ Q_1  (\phi_i+1)^2(\phi_j+1)^2- Q_2 (\phi_i^2-1)^2(\phi_j^2-1)^2\right].
\ee
%where the first term in the function is induced by the repulsive force, the second term is induced by attractive force between cells.

For the   cell-wall adhesion energy, $\phi_w$ is defined along the wall with a  thickness \cite{gu2014simulating} as shown in Figure \ref{fig:init}. Following above interaction potential definition, 
then the cell-wall interaction energy is defined by setting $\phi_2=\phi_w$ in Eq. \eqref{H2phase},
\be
f_w(\phi_i)=  Q_{w_1}  (\phi_i+1)^2(\phi_w+1)^2- Q_{w_2} (\phi_i^2-1)^2(\phi_w^2-1)^2~,
\ee
where $Q_{w_1}$ is repulsive energy density and $Q_{w_2}$ is adhesion energy density.  
Notice that the cell-wall energy  $E_w$ is defined in the bulk not on the boundary.  However, $f_w$ is non-zero only when the two phases overlap which could be regarded as the attraction force is only induced when the vesicle is contacting (or close enough to) the wall. 
% \begin{figure}[ht]
% 	\centering
% 	\includegraphics[width=3.in]{figure/wallattraction/phase.eps}
% 	\caption{A case of the phase definition of wall and vesicle. In the figure, we plot the value of $\phi+\phi_w$. In our model, $\phi$ is 1 and -1 in and outside the vesicle while $\phi_w$ is 0 on the wall and -1 at the edge in the domain. The wall and vesicle phases are defined by tanh function. The thickness of the wall phase is equal to the vesicle membrane. Here, $\phi$ and $\phi_w$ are slightly overlapping which means the attraction force is induced in the overlapping region.}
% 	\label{fig:wall phase}
% \end{figure}

Then the chemical potential $\mu_i$  for each phase $\phi_i$ is defined as follows 
\be\label{mu}
\mu_i = \frac{\delta E_{cell_i}}{\delta \phi_i} = \frac{\hat\kappa_B}{\gamma^3}g(\phi_i)+  \frac {\partial H} {\partial \phi_i}  + \frac{\mathcal{M}_s}{\gamma} \frac {S(\phi_i)-S(\phi_{i,0})} {S(\phi_{i,0})}f(\phi_i)+\frac{\partial f_w(\phi_i)}{\partial \phi_i}~, 
\ee
where $g(\phi_i) = -\gamma^2\Delta f(\phi_i)+(3\phi_i^2-1)f(\phi_i)$.

The dissipation functional of the system consists of the dissipation introduced by fluid viscosity, friction near the wall, and interfacial mixing due to the diffuse interface  representation  \cite{shen2022energy}:
\be\label{dispp}
\Delta &=& \int_{\Omega}2\eta|\bD_{\eta}|^2d\boldsymbol x+\sum_{i=1}^{N}\int_{\Omega}\frac 1 {M_{\phi}}|q_{\phi_i}|^2d\boldsymbol x+  \int_{\partial\Omega_w}\beta_s|\vel_{\btau}|^2ds\nonumber\\
&&+\sum_{i=1}^N\left(\int_{\partial\Omega_w}\kappa_{\Gamma}|J_{\Gamma_i}|^2ds +\int_{\partial\Omega_w}\xi|\gamma \phi_i\nabla\lambda_i|^2ds\right)~, 
\ee
where $\boldsymbol{D}_\eta = \frac{\nabla\vel +(\nabla\vel)^T}{2}$, $\eta$ is the viscosity of fluid mixture,  $\beta_s$ is wall fraction, $M_{\phi}$ and $\kappa_{\Gamma}$ are the mobility of phase in the bulk and boundary. 
We note that in general, the viscosity $\eta$  could be  a function of all phases $\phi_i$ 
\be\label{eta}
\eta = \eta(\phi_w, \phi_1, \cdots, \phi_n)~.
\ee

% Note that the wall phase is taken into consideration since the effect of fraction and block in boundary flow would highly increase under the influence of the wall structure.

The specific forms of the flux and stress functions in the kinematic equations (\ref{assumption})-(\ref{assumption_bd}) are   obtained by taking the time derivative of the total energy functional and comparing with the  defined dissipation functional.   Then the time derivative of the total energy goes:
\be\label{dEdt}
\frac{d E_{total}}{dt} &=& \frac{d}{dt} E_{kin} + \frac{d}{dt} E_{cell} + \frac{d}{dt} E_w\\
& \equiv & I_1+I_2+I_3~.
\ee

Taking the time derivative of $E_{kin}$, together with the conservation law of momentum $(\ref{assumption})_2$, incompressibility of the fluid  $(\ref{assumption})_3$ and 
inextensibility of the membrane $(\ref{assumption})_4$ yields
\be\label{multi_I_1}
I_1 &=&\frac{d}{dt}\int_{\Omega}\frac{\rho|\vel|^2}{2}d\bx \nonumber\\
&=&\int_{\Omega}\frac 1 2 \frac{\partial\rho}{\partial t}|\vel|^2d\bx +\int_{\Omega} \rho\frac{\partial \vel}{\partial t}\cdot\vel d\bx\nonumber\\
&=&\int_{\Omega}\frac 1 2 \frac{\partial\rho}{\partial t}|\vel|^2d\bx +\int_{\Omega} \rho\frac{d \vel}{d t}\cdot\vel d\bx-\int_{\Omega} \left(\rho \vel\cdot\nabla \vel\right)\cdot \vel d\bx\nonumber\\
&=&\int_{\Omega}\frac 1 2 \frac{\partial\rho}{\partial t}|\vel|^2d\bx +\int_{\Omega} \rho\frac{d \vel}{d t}\cdot\vel d\bx+\int_{\Omega} \nabla\cdot(\rho\vel)\frac{|\vel|^2}2d\bx \nonumber\\
&=&\int_{\Omega}( \nabla \cdot \bsigma_\eta)  \cdot\vel d\bx +\int_{\Omega}\bF\cdot \vel d\bx +\sum_{i=1}^{N}\int_\Omega \lambda_i\delta_i \mathcal{P}_i : \na\vel d\bx +\sum_{i=1}^{N}\int_\Omega \xi\gamma^2 \lambda_i\na\cdot(\phi_i^2\na\lambda_i) d\bx-\int_{\Omega} pI:\nabla \vel d\bx \nonumber\\
&=&-\int_{\Omega}((\bsigma_{\eta}+pI):\nabla\vel)d\bx+\int_{\Omega}\bF\cdot \vel dx -\sum_{i=1}^{N}\int_\Omega \na\cdot(\lambda_i\delta_i  \mathcal{P}_i)\cdot\vel d\bx \nonumber\\
&&-\sum_{i=1}^{N}\int_\Omega \xi\gamma^2 \phi_i^2(\na\lambda_i)^2 d\bx+\int_{\partial\Omega_w}((\bsigma_{\eta}+\sum_{i=1}^{N}\lambda_i\delta_i \mathcal{P}_i)\cdot\bn)\cdot\vel_{\tau} dS~,
\ee
where the slip boundary condition is used. Here the pressure $p$ is a Lagrangian multiplier and is introduced to ensure the incompressibility of the fluid.

Taking the time derivative of $E_{cell}$ together with   the conservation of each phase $(\ref{assumption})_1$ yields
\be\label{multi_I_2}
I_2 &=& \frac{d}{dt}\int_{\Omega}\frac{\hat\kappa_B}{2\gamma}\sum_{i=1}^{N}\left|\frac{f(\phi_i)}{\gamma}\right|^2d\bx+\frac{d}{dt}\int_{\Omega}   H d\bx + \sum_{i=1}^{N}\frac {d}{dt}\frac {\mathcal{M}_s} {2} \frac {(S(\phi_i)-S(\phi_{i,0}))^2} {S(\phi_{i,0})}\nonumber \\
&=& \int_{\Omega}\frac{\hat\kappa_B}{\gamma}\sum_{i=1}^{N}\frac{f_i}{\gamma^2}\frac{\partial f_i}{\partial t}d\bx +\int_{\Omega} \sum_{i = 1}^{N}\frac{\partial  H }{\partial \phi_i}\frac{\partial \phi_i}{\partial t}d\bx+\sum_{i=1}^{N} \frac {d}{dt}\frac {\mathcal{M}_s} {2} \frac {(S(\phi_i)-S(\phi_{i,0}))^2} {S(\phi_{i,0})}\nonumber\\
&=& \int_{\Omega}\frac{\hat\kappa_B}{\gamma}\sum_{i=1}^{N}\frac{f_i}{\gamma^2}\left(-\gamma^2\Delta\left(\frac{\partial \phi_i}{\partial t}\right) +(3\phi_i^2-1)\frac{\partial \phi_i}{\partial t}\right)d\bx+\int_{\Omega} \sum_{i = 1}^{N}\frac{\partial  H }{\partial \phi_i}\frac{\partial \phi_i}{\partial t}d\bx\nonumber\\
&&+\sum_{i=1}^{N} \frac {d}{dt}\frac {\mathcal{M}_s} {2} \frac {(S(\phi_i)-S(\phi_{i,0}))^2} {S(\phi_{i,0})}\nonumber\\
&=& \int_{\Omega}\frac{\hat\kappa_B}{\gamma^3}\sum_{i=1}^{N}\left(-\gamma^2\Delta f_i +(3\phi_i^2-1)f_i+ \frac{\partial H}{\partial \phi_i}+ \frac{\mathcal{M}_s}{\gamma} \frac {S(\phi_i)-S(\phi_{i,0})} {S(\phi_{i,0})}f(\phi_i)\right)\frac{\partial \phi_i}{\partial t}d\bx\nonumber\\ 
&&-\int_{\partial\Omega_w}\frac{\hat\kappa_B}{\gamma}\sum_{i=1}^{N}f_i\frac{\partial}{\partial t}(\partial_n\phi_i) ds +\int_{\partial\Omega_w} \frac{\hat\kappa_B}{\gamma} \sum_{i=1}^{N}\partial_n f_i\frac{\partial \phi_i}{\partial t}ds+ \sum_{i=1}^{N}\int_{\partial\Omega_w}\mathcal{M}_s\frac{S(\phi_i)-S(\phi_{i,0})}{S(\phi_{i,0})}\gamma\partial_n \phi_i\frac{\partial \phi_i}{\partial t}ds\nonumber \\
&=&\sum_{i=1}^{N}\left (\int_{\Omega}\mu_i\frac{\partial \phi_i}{\partial t} d\bx +\int_{\partial\Omega_w}\frac{\hat\kappa_B}{\gamma}f_i\frac{\partial}{\partial t}(\partial_n\phi_i) ds+\int_{\partial\Omega_w}\mathcal{M}_s\frac{S(\phi_i)-S(\phi_{i,0})}{S(\phi_{i,0})}\gamma\partial_n \phi_i\frac{\partial \phi_i}{\partial t}ds\right )\nonumber\\
&=&  \sum_{i=1}^{N}\left (\int_{\Omega}\mu_i\nabla \cdot q_{\phi_i}d\bx-\int_{\Omega}\mu_i\vel\cdot\nabla \phi_i d\bx+\int_{\partial\Omega_w} \frac{\hat\kappa_B}{\gamma} \partial_n f_i\frac{\partial \phi_i}{\partial t}ds +\int_{\partial\Omega_w}\mathcal{M}_s\frac{S(\phi_i)-S(\phi_{i,0})}{S(\phi_{i,0})}\gamma\partial_n \phi_i\frac{\partial \phi_i}{\partial t}ds\right )\nonumber\\ 
&=& \sum_{i=1}^{N}\left (-\int_{\Omega}q_{\phi_i}\cdot\nabla \mu_i d\bx-\int_{\Omega}\mu_i\vel\cdot\nabla \phi_i d\bx+\int_{\partial\Omega_w} \frac{\hat\kappa_B}{\gamma} \partial_n f_i\frac{\partial \phi_i}{\partial t}ds\right.\nonumber\\ &&\left.+\int_{\partial\Omega_w}\mathcal{M}_s\frac{S(\phi_i)-S(\phi_{i,0})}{S(\phi_{i,0})}\gamma\partial_n \phi_i\frac{\partial \phi_i}{\partial t}ds\right )~,
\ee 
where the Allan-Cahn boundary condition $(\ref{assumption_bd})_3$ for each phase is used. 

Computing $\frac{d}{dt} E_w$ yields 
\be\label{multi_I_3}
I_3 = \sum_{i=1}^{N}\left (\int_{\partial\Omega}\frac{\partial f_w(\phi_i)}{\partial \phi_i}\frac{\partial \phi_i}{\partial t}ds \right )~. 
\ee

By combining Eqs. \eqref{multi_I_1}, \eqref{multi_I_2} and \eqref{multi_I_3},  we have:
\be
\label{dEdt_expand}
\frac{d}{dt}E_{total}
&=&-\int_{\Omega}((\bsigma_{\eta}+pI):\nabla\vel)d\bx +\int_{\Omega}(\bF-\sum_{i=1}^{N}\mu_i\nabla\phi_i-\sum_{i=1}^{N}\na\cdot(\lambda_i\delta_i \mathcal{P}_i))\cdot \vel d\bx  \nonumber\\
&&- \sum_{i=1}^{N}\int_{\Omega}q_{\phi_i}\cdot \nabla\mu_i d\bx+\sum_{i=1}^{N}\int_\Omega \xi\gamma^2 \phi_i^2(\na\lambda_i)^2 d\bx+\int_{\partial\Omega_w}((\bsigma_{\eta}+\sum_{i=1}^{N}\lambda_i\delta_i  \mathcal{P}_i)\cdot\bn)\cdot\vel_{\tau} ds  +\sum_{i=1}^{N}\int_{\partial\Omega_s}\hat L_i\frac{\partial \phi_i}{\partial t} ds \nonumber\\
&=&-\int_{\Omega}((\bsigma_{\eta}+pI):\nabla\vel)d\bx+\int_{\Omega}(\bF-\sum_{i=1}^{N}\mu_i\nabla\phi_i-\sum_{i=1}^{N}\na\cdot(\lambda_i\delta_i \mathcal{P}_i))\cdot \vel d\bx \nonumber\\
&&- \sum_{i=1}^{N}\int_{\Omega}q_{\phi_i}\cdot \nabla\mu_i d\bx+\sum_{i=1}^{N}\int_\Omega \xi\gamma^2 \phi_i^2(\na\lambda_i)^2 d\bx+\int_{\partial\Omega_w}((\bsigma_{\eta}+\sum_{i=1}^{N}\lambda_i\delta_i  \mathcal{P}_i)\cdot\bn)\cdot\vel_{\tau} ds\nonumber\\
&&+\sum_{i=1}^{N}\int_{\partial\Omega_s}\hat L_i(-\vel\cdot\nabla_{\Gamma}\phi_i+J_{\Gamma_i})ds\nonumber\\
&=&-\int_{\Omega}((\bsigma_{\eta}+pI):\nabla\vel)d\bx+\int_{\Omega}(\bF-\sum_{i=1}^{N}\mu_i\nabla\phi_i-\sum_{i=1}^{N}\na\cdot(\lambda_i\delta_i \mathcal{P}_i))\cdot \vel d\bx \nonumber\\
&&- \sum_{i=1}^{N}\int_{\Omega}q_{\phi_i}\cdot \nabla\mu_i d\bx+\sum_{i=1}^{N}\int_\Omega \xi\gamma^2 \phi_i^2(\na\lambda_i)^2 d\bx+\int_{\partial\Omega_w}((\bsigma_{\eta}+\sum_{i=1}^{N}\lambda_i\delta_i  \mathcal{P}_i)\cdot\bn-\sum_{i=1}^{N}\hat L_i\nabla_{\Gamma}\phi_i)\cdot\vel_{\tau} ds \nonumber\\
&&+\sum_{i=1}^{N}\int_{\partial\Omega_w}\hat L_i J_{\Gamma_i}ds~,
\ee
where $ \hat L_i = \displaystyle\frac{\hat\kappa_B}{\gamma}\partial_n f_i +\mathcal{M}_s\frac {S(\phi_i)-S(\phi_{i,0})}{S(\phi_{i,0})}\gamma \partial_n\phi_i$.

To close the system, we use the energy dissipation law \cite{xu2014energetic,shen2020energy, liu2019energetic}
\be
\frac{dE}{dt} = -\Delta~,
\ee
which states that the rate of changing total energy of the system is induced by the dissipation.
Comparing Eq.~(\ref{dEdt_expand}) with predefined dissipation functional in  Eq. \eqref{dispp} gives rise to
\be\left\{\begin{array}{ll}
	\bsigma_{\eta} = 2\eta\bD_{\eta}-pI~,& \mbox{in~} \Omega~,\\
	q_{\phi_j} = M_{\phi_j}\nabla\mu_j~, &\mbox{in~} \Omega~,\\
	\bF =\sum_{i=1}^{N}( \mu_i\nabla\phi_i+\nabla\cdot(\lambda_i\delta_i \mathcal{P}_i))~, & \mbox{in~}\Omega~,\\
	J_{\Gamma_i} =- \kappa_{\Gamma_i}^{-1} \hat L_i~, &\mbox {on~}\partial\Omega_w~,\\
	u_{\tau_j} =\beta_{s}^{-1}(-(\bn\cdot(\bsigma_{\eta}+\sum_{i=1}^{N}\lambda_i\delta_i \mathcal{P}_i)\cdot\btau_j)+\sum_{i=1}^{N}\hat L_i\partial_{\tau_j}\phi_i)~,~ j=1,2,& \mbox{on~} \partial\Omega_w~.
\end{array}\right.
\ee

To this end, the proposed multi-cellular interaction model is composed of the following
equations:

\be\left\{\begin{array}{lll}
	\label{assumption_det}
	\frac{\partial \phi_i}{\partial t} +\nabla\cdot (\vel \phi_i) = M_{\phi_i}\Delta \mu_i~,\\
	\mu_i = \frac{\hat\kappa_B}{\gamma^3}g(\phi_i)+ \frac{\partial  H}{\partial \phi_i}+\frac{M_s}{\gamma}\frac{S(\phi_i)-S(\phi_{i,0})}{S(\phi_{i,0})}f(\phi_i) +\frac{\partial f_w(\phi_i)}{\partial \phi_i}~,\\
	g(\phi_i) = -\gamma^2\Delta f_i+(3\phi^2-1)f(\phi_i),\\ f(\phi_i)= -\gamma^2 \Delta\phi_i+ (\phi_i^2-1)\phi_i~,\\
	\rho(\frac{\partial\vel}{\partial t}+(\vel\cdot\nabla)\vel)+\nabla p = \nabla\cdot(2\eta\bD_{\eta})+\sum_{j=1}^{i}( \mu_j\nabla\phi_j+\nabla\cdot(\lambda_j\delta_j \mathcal{P}_j))~,\\
	\nabla \cdot\vel = 0~,\\
	\delta_i ( \mathcal{P}_i : \na\vel)+\xi\gamma^2 \nabla\cdot(\phi_i^2\nabla\lambda_i)=0~,
\end{array}\right.
\ee
with the boundary conditions on $\partial\Omega_w$
\be\left\{\begin{array}{l}
	\label{sys_bd}
	\vel\cdot\bn =0~,\\
	-\beta_{s}u_{\tau_j} =(\bn\cdot(\bsigma_{\eta}+\sum_{i=1}^{N}\lambda_i\delta_i \mathcal{P}_i)\cdot\btau_j)-\sum_{i=1}^{N}\hat L_i\partial_{\tau_j}\phi_i~, ~j=1,2,\\
	f_i=0~,\\
	\kappa_{\Gamma_i}\left(\frac{\partial \phi_i}{\partial t} + \vel\cdot\nabla_{\Gamma}\phi_i\right) = -\hat L_i~,\\
	\hat L_i = \frac{\hat\kappa_B}{\gamma}\partial_n f_i +\mathcal{M}_s\frac {s(\phi_i)-S(\phi_{i,0})}{S(\phi_{i,0})}\gamma \partial_n\phi_i~,\\
	\partial_n\mu_i = 0~.
\end{array}
\right.
\ee

\subsection{Dimensionless governing equations}
The  viscosity, length, velocity, time, bulk and boundary chemical potentials in the equations are scaled by their corresponding characteristic values $\eta_0$, $L$, $U$, $\frac{L}{U}$, $\frac{\eta_0 U}{L}$ and $\eta_0U$, respectively.  
Write $Q_{w_1}, Q_{w_2}, Q_1, Q_2$ into $Q_0q_{w_1}, Q_0q_{w_2}. Q_0q_1, Q_0q_2$, where $Q_0$ is the character energy density. 
The governing equation of the system can be rewritten as 

\be\left\{\begin{array}{ll}
	\label{sys_nd}
	Re(\frac{\partial \vel}{\partial t}+(\vel\cdot \na)\vel)+\na P=\na\cdot(2\eta \bD)+\sum_i\mu_i\na\phi_i+\sum_i\na\cdot(\lambda_i\delta_{\epsilon_i} \mathcal{P}_i)~,& \mbox{in~}\Omega~,\\[3mm]
	\na\cdot \vel=0~,& \mbox{in~}\Omega~,\\[3mm]
	\frac{\partial\phi_i}{\partial t}+\vel\cdot\na\phi_i=-\mathcal{M}\Delta\mu_i~,& \mbox{in~}\Omega~,\\[3mm]
	\mu_i=     \kappa_B g(\phi_i) +\mathcal{M}_s\frac{(S(\phi_i)-S(\phi_{i,0}))}{S(\phi_{i,0})}f_i+\alpha \frac{\partial H}{\partial \phi_i}+\alpha\frac{\partial f_w(\phi_i)}{\partial\phi_i}~,  & \mbox{in~}\Omega~,\\[3mm]
	f_i = -\epsilon \Delta\phi_i+ \frac{(\phi_i^2-1)}{\epsilon}\phi_i,~g(\phi_i) = -\Delta f_i+\frac{1}{\epsilon^2}(3\phi_i^2-1)f_i~, & \mbox{in~}\Omega~,\\[3mm]
	\delta_{\epsilon_i}( \mathcal{P}_i : \na\vel)+\xi\eps^2 \nabla\cdot(\phi_i^2\nabla\lambda_i)=0 ~, & \mbox{in~}\Omega~,
\end{array}\right.
\ee
with  the boundary conditions 
\be
\displaystyle \left\{\begin{array}{ll}
	\label{bd_nd}
	\kappa \dot{\phi_i}+L(\phi_i)=0~, & \mbox{on~} \partial\Omega_w~,\\
	L(\phi_i)= \kappa_B \partial_n f(\phi_i)+ \epsilon\mathcal{M}_s\frac{S(\phi_i)-S(\phi_{i,0})}{S(\phi_{i,0})} \partial_n\phi_i~,  & \mbox{on~} \partial\Omega_w~,\\
	-l_s^{-1} u_{\tau_i} = \boldsymbol{\tau_i}\cdot( 2\eta \mathbf{D}_{\eta}+\sum_i\lambda_i\delta_{\epsilon_i} \mathcal{P}_i)\cdot \boldsymbol{n} - \sum_i L(\phi_i)\partial_{\tau_i}\phi_i~,~i=1,2, & \mbox{on~}  \partial\Omega_w~,\\
	f_i=0~, &\mbox{on~} \partial\Omega_w~,\\
	\partial_n\mu_i = 0~, &\mbox{on~} \partial\Omega_w~,
\end{array}\right.
\ee

where $S(\phi_i) = \displaystyle\int_{\Omega}\frac{\epsilon}{2}|\nabla\phi_i|^2 +\frac{1}{4\epsilon}(\phi_i^2-1)^2 d\bx$ and $ \delta_{\epsilon_i}=\frac{1}{2}\epsilon^2|\na\phi_i|^2 $.

The dimensionless constants appeared in Eqs.~(\ref{sys_nd})-(\ref{bd_nd}) are given by $\epsilon = \frac{\gamma}{L}$, $Re = \frac{\rho_0 UL}{\eta_0}$, $\mathcal{M} =\frac {M_\phi\eta_0}{L^2}$, 
$\kappa_B = \frac{\hat\kappa_B}{L^2\eta_0U}$, 
$k = \frac{\hat\kappa_B  }{\eta_0 L}$, 
$l_s = \frac{\eta_0}{\beta_{s}L}$,  $\alpha=\frac{Q_0}{\eta_0 U}$, $\mathcal{M}_s = \frac{M_s}{\eta_0 U}$.  

If we define the Sobolev spaces as follows 
\cite{guo2014numerical,shen2020energy}
\begin{eqnarray}
\boldsymbol{W}^{1,3}=(W^{1,3})^2~,\\
\boldsymbol{W}_N^{1,\frac 3 2}=(W^{1,\frac 3 2})^N~,\\
\boldsymbol{W}_{N}^{1,3}=(W^{1,3})^N~,\\
\boldsymbol{W}_N^{1,\frac 3 2}(\Omega)=\left\{\Lambda=(\lambda_1,\lambda_2,...\lambda_N)^T \right\}~,\\
\boldsymbol{W}^{1,3}(\Omega) = \left\{\vel=(u_x,u_y)^T\in  \boldsymbol{W}^{1,3}|\vel\cdot\boldsymbol{n}=0,~\mbox{on~} \partial\Omega_w \right\}~, \\
\boldsymbol{W}_\Phi^{1,3}(\Omega) = \left\{\Phi=(\phi_1,\phi_2,...\phi_N)^T\in  \boldsymbol{W}_N^{1,3}|-1\leq \phi_i\leq  1,i=1,2,...N,~\mbox{in~}\Omega\right\}~, \\
\boldsymbol{W}_{\boldsymbol{U}}^{1,3}(\Omega) = \left\{\boldsymbol{U}=(\mu_1,\mu_2,...\mu_N)^T\in  \boldsymbol{W}_N^{1,3}|\partial_n \mu_i=0,i=1,2,...N,~\mbox{on~}\partial\Omega_w\right\}~, \\
\boldsymbol{W}_F^{1,3}(\Omega) = \left\{F=(f_1,f_2,...f_N)^T\in  \boldsymbol{W}_N^{1,3}| f_i=0,i=1,2,...N,~\mbox{on~}\partial\Omega_w\right\}~, \\
\boldsymbol{W}_b =  \boldsymbol{W}_\Phi^{1,3}(\Omega)\times\boldsymbol{W}_F^{1,3}(\Omega)\times \boldsymbol{W}_{\boldsymbol{U}}^{1,3}(\Omega)\times   	\boldsymbol{W}_N^{1,\frac 3 2}(\Omega)\times W^{1,\frac 3  2}(\Omega)  \times\boldsymbol{W}^{1,3}(\Omega) ~,
\end{eqnarray}
and let  $\| \cdot \| = \left(\int_{\Omega}| \cdot |^2 d\bx\right)^{\frac 1 2} $  and $\| \cdot \|_w = \left(\int_{\partial\Omega_w}| \cdot |^2 ds\right)^{\frac 1 2} $ denote the $L^2$ norm defined in the domain and on the domain boundary respectively,  then the system \eqref{sys_nd}-\eqref{bd_nd} satisfies the following energy law.

\begin{theorem}\label{thm:edp}
	If $(\Phi,~F,~\boldsymbol{U},~\lambda, ~P,~\vel)\in \boldsymbol{W}_b$ are  smooth solutions of the above system \eqref{sys_nd}-\eqref{bd_nd}, then  the following energy law is satisfied:
	\be\label{energye3}
	&&\frac{d}{dt}\mathcal{E}_{total}=\frac{d}{dt}(\mathcal{E}_{kin}+\mathcal{E}_{cell}+\mathcal{E}_w)\nonumber\\
	&=&\frac{1}{Re}\left( -2\|\eta^{1/2}\mathbf{D}_{\eta}\|^2  - \mathcal{M}\sum_i\|\na\mu_i\|^2-\xi\sum_i\|\epsilon\phi_i\nabla\lambda_i\|^2- \kappa \sum_i\|\dot{\phi_i}\|_{w}^2-\|l_s^{-1/2} \vel_{\tau}\|^2_{w}\right)~,
	\ee
	where 
	$\mathcal{E}_{total}=\mathcal{E}_{kin}+\mathcal{E}_{cell}+\mathcal{E}_{w}$,
	$\mathcal{E}_{kin}= \displaystyle\frac{1} 2 \int_{\Omega}|\vel|^2 d\bx $,
	$\mathcal{E}_{cell} =  \displaystyle\frac{\kappa_B }{2Re\epsilon} \sum_i\int_{\Omega}|f_i|^2 d\bx+\mathcal{M}_s\sum_i\frac{(S(\phi_i)-S(\phi_{i,0}))^2}{2Re S(\phi_{i,0})}+\frac{\alpha}{Re}\int_{\Omega}H d\bx$
	and
	$\mathcal{E}_w\!\! =\!\! \frac{\alpha}{Re}\!\!  \sum_i\displaystyle\int_{\Omega}f_{w}(\phi_i) d\bx$. 
\end{theorem}

\noindent \textbf{Proof:} Multiplying the first equation in Eq.~\eqref{sys_nd}  with $\vel$ and integration by parts yield 
\be\label{kin_nd}
\frac{d}{dt}\mathcal{E}_{kin}&=&\frac{1}{Re}\left\{-\int_{\Omega}2\eta|\bD_{\eta}|^2d\bx+\int_{\partial\Omega_w}(\bsigma_{\eta}\cdot\bn)\cdot \vel_{\tau}ds +\sum_i\int_{\Omega}\mu_i\nabla\phi_i\cdot\vel d\bx-\sum_i\int_{\Omega}\lambda_i\delta_{\epsilon_i} \mathcal{P}_i : \nabla \vel d\bx\right.\nonumber\\
&&\left.+\sum_i\int_{\partial\Omega_w}( \lambda_i\delta_{\epsilon_i} \mathcal{P}_i\cdot \bn) \cdot\vel_{\tau} ds\nonumber\right\}\\
&=&\frac{1}{Re}\left\{-\int_{\Omega}2\eta|\bD_{\eta}|^2d\bx-\sum_i\int_{\Omega}\lambda_i\delta_{\epsilon_i} \mathcal{P}_i : \nabla \vel d\bx-l_s^{-1}\int_{\partial\Omega_w}|\vel_{\tau}|^2ds\right.\nonumber\\
&&\left.+\sum_i\int_{\partial\Omega_w} L(\phi_i)\partial_{\tau}\phi\cdot\vel_{\tau}ds+\sum_i\int_{\Omega}\mu_i\nabla\phi_i\cdot\vel d\bx~,\right\}
\ee 
where  the slip boundary condition in Eq.~\eqref{bd_nd} is applied.

%\textcolor{black}{(I am lost in the grammar here.)}

Taking the inner product of the third equation in Eq.~\eqref{sys_nd}   with $\frac{\mu_i}{Re}$ and summing up with respect to $i$ result in

\be\label{phit1}
\frac{1}{Re}\sum_i\int_{\Omega}\frac{\partial\phi_i}{\partial t}\mu_i d\bx +\frac{1}{Re}\sum_i\int_{\Omega}\vel\cdot\nabla\phi_i\mu_i d\bx =-\frac{1}{Re}\mathcal{M}\sum_i\int_{\Omega}|\nabla\mu_i|^2d\bx~,
\ee
where $\partial_n \mu_i=0$ is considered here. 

Multiplying the fourth equation in Eq.~\eqref{sys_nd}  with $\frac 1 {Re}\frac{\partial\phi_i}{\partial t}$ and integration by parts give rise to 
\be\label{mut}
&&\quad \frac 1 {Re}\sum_i\int_{\Omega}\mu\frac{\partial\phi_i}{\partial t}d\bx\\
&=&\frac 1 {Re}\sum_i\left\{\kappa_B\int_{\Omega}g_i\frac{\partial\phi_i}{\partial t}d\bx +\mathcal{M}_s\frac{S(\phi_i)-S(\phi_{i,0})}{S(\phi_{i,0})}\int_{\Omega}f_i\frac{\partial\phi_i}{\partial t}d\bx +\alpha\int_{\Omega}\frac{\partial H}{\partial \phi_i}\frac{\partial\phi_i}{\partial t}d\bx+\alpha\int_{\Omega}\frac{f_w(\phi_i)}{\partial\phi_i}\frac{\partial\phi_i}{\partial t}d\bx\right\}\nonumber\\
&=&\frac {\kappa_B} {Re}\sum_i\int_{\Omega}f_i\frac{\partial}{\partial t}\left(-\Delta \phi_i+\frac{1}{\epsilon^2}(\phi_i^3-\phi_i)\right)d\bx-\frac {\kappa_B} {Re}\sum_i\int_{\partial\Omega_w}\partial_n f_i\frac{\partial\phi_i}{\partial t}ds+\mathcal{M}_s\sum_i\frac{d}{dt}\left(\frac{(S(\phi_i)-S(\phi_{i,0}))^2}{2Re S(\phi_{i,0})}\right)\nonumber\\
&&-\mathcal{M}_s\sum_i\left(\frac{S(\phi_i)-S(\phi_{i,0})}{ Re S(\phi_{i,0})}\right)\int_{\partial\Omega_w}\epsilon\partial_n\phi_i\frac{\partial\phi_i}{\partial t} ds+\frac{\alpha}{Re}\sum_i\int_{\Omega}\frac{\partial f_w(\phi_i)}{\partial\phi_i}\frac{\partial\phi_i}{\partial t}d\bx +\frac{\alpha}{Re}\sum_i\int_{\Omega}\frac{\partial H}{\partial \phi_i}\frac{\partial\phi_i}{\partial t}d\bx \nonumber\\
&=&\frac{d}{dt}\left(\kappa_B\sum_i\int_{\Omega}\frac{|f_i|^2}{2Re \epsilon}d\bx\right)+\mathcal{M}_s\frac{d}{dt}\left(\sum_i\frac{(S(\phi_i)-S(\phi_{i,0}))^2}{2 Re S(\phi_{i,0})}\right)+\frac{\alpha}{Re}\frac{d}{dt}\sum_i\int_{\Omega} f_w(\phi_i) d\bx+\frac{\alpha}{Re}\frac{d}{dt}\int_{\Omega} H d\bx\nonumber\\
&&-\sum_i\int_{\partial\Omega_w}\frac{ L(\phi_i)}{Re} \frac{\partial\phi_i}{\partial t}ds \nonumber\\
&=&\frac{d}{dt}(\mathcal{E}_{cell}+ \mathcal{E}_w )-\int_{\partial\Omega_w} \frac{L(\phi)}{Re} \frac{\partial\phi}{\partial t}ds\nonumber~,
\ee
where  the definitions of $f(\phi)$, $g(\phi)$ and the boundary conditions of $\phi$ and $f$ are utilized. 

Multiplying the last equations with $ \frac{\lambda_i}{Re}$  and integration by parts and sum up by $i$ leads to
\be\label{lamda}
&&   \frac{1}{Re} \sum_i\int_\Omega (\lambda_i\delta_{\epsilon_i} \mathcal{P}_i)  :\nabla\vel d\bx-\frac{1}{Re} \sum_i\int_\Omega \xi\eps^2\phi_i^2(\na\lambda_i)^2 d\bx
%-  \frac{1}{Re} \int_\Omega \xi\eps^2\phi^2(\na\lambda)^2%
=0~.
\ee

Finally, the energy dissipation law \eqref{energye3} is obtained by combining Eqs. \eqref{kin_nd}, \eqref{phit1}, \eqref{mut} and \eqref{lamda} considering the boundary conditions in \eqref{bd_nd}. 
$\null \hfill \blacksquare$

\section{Numerical Scheme and Discrete Energy Law}
\label{sec:numericalscheme}

In this section, an energy stable temporal discretization scheme is first proposed for the multi-cellular system \eqref{sys_nd}-\eqref{bd_nd}. Then a $C^0$ finite element method is used for spacial discretization to obtain the fully discrete scheme.

\subsection{Time-discrete primitive method}
The mid-point method is used  for the temporal discretization of Eqs.~\eqref{sys_nd}-\eqref{bd_nd}. Let $\Delta t$ denote the time step size, $()^{n+1}$ and $()^{n}$ denote the values of the variables at times $(n+1)\Delta t$ and $n\Delta t$, respectively. The semi-discrete  in time scheme to solve Eq.~(\ref{sys_nd}) is as follows:
%in $\Omega$
\be
\left\{\begin{array}{ll}\label{sys_nd_dis}
	\frac{\vel^{n+1}-\vel^{n}}{\Delta t}+(\vel^{n+\frac{1}{2}}\cdot\na) \vel^{n+\frac{1}{2}}+\frac{1}{Re}\na P^{n+\frac{1}{2}}=\frac{1}{Re}\na\cdot(\eta^{n}(\na \vel^{n+\frac{1}{2}} +(\na \vel^{n+\frac{1}{2}})^T))\\ 
	~~~~~~~~	+\frac{1}{Re} \sum_{i} \mu_i^{n+\frac{1}{2}}\na\phi_i^{n+\frac{1}{2}}+\sum_i\frac{1}{Re}\na\cdot \left(\lambda_i^{n+\frac{1}{2}} \mathcal{P}_i^{n} \delta_{\epsilon_i}\right)~, \\ 
	\na\cdot \vel^{n+\frac{1}{2}}=0~,\\ 
	\frac{\phi_i^{n+1}-\phi_i^{n}}{\Delta t}+(\vel^{n+\frac{1}{2}}\cdot\na)\phi_i^{n+\frac{1}{2}}=-\mathcal{M}\Delta\mu_i^{n+\frac{1}{2}}~,\\ 
	\mu_i^{n+\frac{1}{2}} =     \kappa_B g(\phi_i^{n+1},\phi_i^n)
	+\mathcal{M}_s\frac{(S(\phi_i^{n+\frac{1}{2}})-S(\phi_{0_i}))}{S(\phi_{0_i})}f(\phi_i^{n+1},\phi_i^n)\\
	~~~~~~~~~	+\alpha\frac{ H_i^{n+1}-H_i^n}{\phi_i^{n+1  }-\phi_i^{n }}+\alpha\frac{f_w(\phi_i^{n+1})-f_w(\phi_i^n)}{\phi_i^{n+1}-\phi_i^n}~, \\ 
	f_i^{n+\frac{1}{2}} = -\epsilon \Delta\phi_i^{n+\frac{1}{2}}+\frac{1}{\epsilon}((\phi_i^{n+\frac{1}{2}})^2-1)\phi_i^{n+\frac{1}{2}}~, \\ 
	%\xi \epsilon^2 \nabla \cdot((\phi_i^{n})^2\nabla \lambda_i^{n+\frac1 2})+	\delta_{\epsilon_i} \mathcal{P}_i^{n} : \na\vel^{n+\frac{1}{2}}=0~.
	\delta_{\epsilon_i} \mathcal{P}_i^{n} : \na\vel^{n+\frac{1}{2}}+\xi \epsilon^2 \nabla \cdot((\phi_i^{n})^2\nabla \lambda_i^{n+\frac1 2})=0~,
\end{array}\right.
\ee
with boundary conditions on $ \partial\Omega_w$~,
\be\left\{\begin{array}{ll}\label{bd_nd_dis}
	\kappa\dot \phi_i^{n+\frac{1}{2}}=-L_i^{n+\frac{1}{2}}~, \\
	L_i^{n+\frac{1}{2}}= \kappa_B \partial_n f_i^{n+\frac{1}{2}}+ \mathcal{M}_s \epsilon\frac{S(\phi_i^{n+\frac{1}{2}})-S_{0_i}}{S_{0_i}}\partial_n\phi_i^{n+\frac{1}{2}}~,  \\
	-l_s^{-1} u_{\tau_j}^{n+\frac{1}{2}} = \boldsymbol{\tau_j}\cdot (\eta^{n}(\na \vel^{n+\frac{1}{2}} +(\na \vel^{n+\frac{1}{2}})^T)+\sum_i\lambda_i^{n+\frac{1}{2}}\delta_{\epsilon_i} \mathcal{P}_i^{n})\cdot \boldsymbol{n} \\
	~~~~~~~~~~~~~~~~ - \sum_i L_i^{n+\frac{1}{2}}\partial_{\tau_j}\phi_i^{n+\frac{1}{2}} ~,~~j=1,2, \\
	f_i^{n+\frac{1}{2}}=0~, \\
	\partial_n \lambda_i^{n+\frac 1 2} =0~,
\end{array}\right.\ee
with
$(\cdot)^{n+\frac 1 2} = \frac {(\cdot)^n +(\cdot)^{n+1}}{2}$ and $\mathcal {P}_i^n = I -\bn^{n}_m\otimes\bn^{n}_m$ with $\bn_m^n = \frac{\nabla\phi_i^{n}}{|\nabla\phi_i^n|}$ and
\be\left\{\begin{array}{l}
	f(\phi_i^{n+1},\phi_i^n) = -\epsilon \Delta \phi_i^{n+\frac{1}{2}}+\frac{1}{4\epsilon}((\phi_i^{n+1})^2+(\phi_i^{n})^2-2)(\phi_i^{n+1}+\phi_i^n)~,\\
	g(\phi_i^{n+1},\phi_i^n) = \left(-\Delta f_i^{n+\frac{1}{2}}+\frac{1}{\epsilon^2}\left((\phi_i^{n+1})^2+(\phi_i^{n})^2+\phi_i^{n+1}\phi_i^{n}-1\right)f_i^{n+\frac{1}{2}}\right)~, \\
	H_i^n = q_1  (\phi_i^{n}+1)^2\sum_{j\neq i}\left [(\phi_j^{n+\frac 1 2}+1)^2\right ]-q_2( (\phi_i^{n})^2-1)^2\sum_{j\neq i}\left[((\phi_j^{n+\frac 1 2})^2-1)^2\right]~,\\
	H_i^{n+1} = q_1 (\phi_i^{n+1}+1)^2\sum_{j\neq i}\left [(\phi_j^{n+\frac 1 2}+1)^2\right ]-q_2( (\phi_i^{n+1})^2-1)^2\sum_{j\neq i}\left[((\phi_j^{n+\frac 1 2})^2-1)^2\right]~,\\
	f_w(\phi^n) =q_{w_1}(\phi_i^n+1)^2(\phi_w+1)^2-q_{w_2}((\phi_i^n)^2-1)((\phi_w)^2-1) ~.
\end{array}\right. 
\ee

Thus we have
\be
\frac{H_i^{n+1 }-H_i^n}{\phi_i^{n+1  }-\phi_i^{n }} &=&\frac{1}{\phi_i^{n+1  }-\phi_i^{n }}  \left( {q_1}  (\phi_i^{n+1}+1)^2\sum_{j\neq i}\left [(\phi_j^{n+\frac 1 2}+1)^2\right ]-q_2( (\phi_i^{n+1})^2-1)^2\sum_{j\neq i}\left[((\phi_j^{n+\frac 1 2})^2-1)^2\right]\right.\nonumber\\
&&-\left. q_1 (\phi_i^n+1)^2\sum_{j\neq i}\left [(\phi_j^{n+\frac 1 2}+1)^2\right ]+q_2((\phi_i^n)^2-1)^2\sum_{j\neq i}\left[(\phi_j^{n+\frac 1 2})^2-1\right] \right) \nonumber\\
&=& q_1 (\phi_i^{n+1}+\phi_i^n+2)\sum_{j\neq i}\left [(\phi_j^{n+\frac 1 2}+1)^2\right ]\nonumber\\
&&-q_2 (\phi_i^{n+1}+\phi_i^n)((\phi_i^{n+1})^2+(\phi_i^n)^2-2)\sum_{j\neq i}\left[((\phi_j^{n+\frac 1 2})^2-1)^2\right]
\ee
Similarly,
\be
\frac{f_w(\phi_i^{n+1})-f_w(\phi_i^n)}{\phi_i^{n+1}-\phi_i^n} &=&q_{w_1} (\phi_i^{n+1}+\phi_i^n+2)(\phi_w+1)^2\nonumber\\
&&-q_{w_2} (\phi_i^{n+1}+\phi_i^n)((\phi_i^{n+1})^2+(\phi_i^n)^2-2)(\phi_w^2-1)^2
\ee
Later on, we keep the form $\frac{f_w(\phi_i^{n+1})-f_w(\phi_i^n)}{\phi_i^{n+1}-\phi_i^n}$ and $\frac{H(\phi_i^{n+1})-H(\phi_i^n)}{\phi_i^{n+1}-\phi_i^n}$for convenience in later derivation.

The above scheme obeys the following  theorem of energy stability. 
\begin{theorem}
	\label{energyTh_dis}
	If $(\phi_i^n, \mu_i^n, \vel^n, P^n)$ are  smooth solutions of the above system \eqref{sys_nd_dis}-\eqref{bd_nd_dis}, then  the following energy law is satisfied:
	\be\label{energye3_dis}
	\mathcal{E}_{total}^{n+1}-\mathcal{E}_{total}^n&\!\!\!=\!\!\!& (\mathcal{E}_{kin}^{n+1}+\sum_i^N\left[\mathcal{E}_{cell_i}^{n+1}+\mathcal{E}_{i,int}^{n+1}+\mathcal{E}_{w_i}^{n+1}\right])- (\mathcal{E}_{kin}^{n}+\sum_i^N\left[\mathcal{E}_{cell_i}^{n}+\mathcal{E}_{i,int}^{n}+\mathcal{E}_{w_i}^{n}\right])\nonumber\\
	&=& \frac{ \triangle t}{Re}\left(-2\|(\eta^n)^{1/2}\mathbf{D}^{n+\frac 12}_{\eta}\|^2  - \mathcal{M}\sum_i^N\|\na\mu_i^{n+\frac 1 2}\|^2- \xi\sum_i^N\|\ \epsilon \phi_i^n \nabla \lambda_i^{n+\frac 1 2}  \|^2\right.\nonumber\\%- \xi\|\ \epsilon \phi^n \nabla \lambda^{n+\frac 1 2}  \|^2\right.\nonumber\\
	&&\left.- \frac 1 {\kappa} \|\sum_i^N L(\phi_i^{n+\frac 1 2})\|_{w}^2-\|l_s^{-1/2} \vel_{\tau}^{n+\frac 1 2}\|^2_{w}\right)~,
	\ee
	where 
	$\mathcal{E}_{total}^n=\mathcal{E}_{kin}^n+\sum_i^N\left[\mathcal{E}_{cell_i}^n+\mathcal{E}_{i,int}^n+\mathcal{E}_{w_i}^n\right]$ with
	$\mathcal{E}_{kin}^n=\frac{1} 2 \|\vel^n\|^2$,
	$\mathcal{E}_{cell_i}^n = \frac{\kappa_B\|f_i^n\|^2}{2Re\epsilon}+\mathcal{M}_s\frac{(S(\phi_i^ n)-S(\phi_{i,0}))^2}{2ReS(\phi_{i,0})}+\frac{\alpha}{Re}H_i^n$
	and
	$\mathcal{E}_{i,w}^n\!\! =\!\! \frac{\alpha}{Re}\!\!  \displaystyle\int_{\Omega}f_{i,w}^n d\bx$. 
\end{theorem}
The following two   lemmas are needed for proving \textbf{Theorem \ref{energyTh_dis}}. Proof of these two lemmas can be found in \cite{shen2022energy}.
\begin{lemma}\label{lma:1}
	Let
	\be 
	f(\phi^{n+1},\phi^{n})=-\epsilon \Delta \phi^{n+\frac{1}{2}}+\frac{1}{4\epsilon}((\phi^{n+1})^2+(\phi^{n})^2-2)(\phi^{n+1}+\phi^n)~.
	\ee
	Then $f(\phi^{n+1},\phi^{n})$ satisfies
	\be
	\int_\Omega f(\phi^{n+1},\phi^{n})(\phi^{n+1}-\phi^{n})d\bx=S^{n+1}-S^n-\int_{\partial\Omega_w} \epsilon\partial_n \phi^{n+\frac{1}{2}}(\phi^{n+1}-\phi^{n})ds~,
	\ee
	where $S^{n+1}=\int_\Omega G(\phi^{n+1})d\bx, S^{n}=\int_\Omega G(\phi^{n})d\bx$.
\end{lemma}

\begin{lemma}\label{lma:2}
	Let
	$g(\phi^{n+1},\phi^{n})=-\Delta f^{n+\frac{1}{2}}+\frac{1}{\epsilon^2}((\phi^{n+1})^2+(\phi^{n})^2+\phi^{n+1}\phi^n-1)f^{n+\frac{1}{2}}~.$
	Then $g(\phi^{n+1},\phi^{n})$ satisfies
	\be
	&&\int_\Omega g(\phi^{n+1},\phi^n)(\phi^{n+1}-\phi^n)d\bx\nonumber\\
	&=&\int_\Omega \frac{1}{2\epsilon}((f^{n+1})^2-(f^n)^2)d\bx-\int_{\partial\Omega_w}\partial_n f^{n+\frac{1}{2}}(\phi^{n+1}-\phi^n)ds~,
	\ee
	where $f^{n+1}=-\epsilon \Delta \phi^{n+1}+\frac{1}{\epsilon}((\phi^{n+1})^2-1)\phi^{n+1}, f^{n}=-\epsilon \Delta \phi^{n}+\frac{1}{\epsilon}((\phi^{n})^2-1)\phi^{n}$.
\end{lemma}	

\noindent{\textbf{Proof of Theorem \ref{energyTh_dis}}:}  Multiplying the first equation in system \eqref{sys_nd_dis} by $\Delta t \vel^{n+\frac{1}{2}}$ gives
\be\label{kenetic_dis}
&&\int_{\Omega} \frac{1}{2}((\vel^{n+1})^2-(\vel^{n})^2)d\bx+\int_{\Omega}\Delta t \vel^{n+\frac{1}{2}}\cdot ((\vel^{n+\frac{1}{2}}\na)\cdot \vel^{n+\frac{1}{2}})d\bx\nonumber\\
&&-\frac{\Delta t}{Re}\int_{\Omega}P^{n+\frac{1}{2}}\na\cdot \vel^{n+\frac{1}{2}}d\bx\nonumber\\
&=&-\frac{\Delta t}{Re}\int_{\Omega}\na \vel^{n+\frac{1}{2}}:\eta^{n}(\na \vel^{n+\frac{1}{2}}+(\na \vel^{n+\frac{1}{2}})^T)d\bx+\frac{\Delta t}{Re}\sum_i\int_{\Omega}\vel^{n+\frac{1}{2}}\cdot\na \phi_i^{n+1}\mu_i ^{n+1}d\bx\nonumber\\
&&- \frac{\Delta t}{Re} \sum_i\int_{\Omega}\lambda_i\delta_{\epsilon_i} \mathcal{P}_i^n : \nabla \vel^{n+\frac 1 2} d\bx+\frac{\Delta t}{Re}\sum_i\int_{\partial\Omega_w} \lambda_i^{n+\frac{1}{2}}( \delta_{\epsilon_i}\mathcal{P}_i^{n}\cdot\bn)\cdot \vel_\tau ^{n+\frac{1}{2}} ds\nonumber\\
&&+\frac{\Delta t}{Re}\int_{\partial\Omega_w}\vel^{n+\frac{1}{2}}\cdot\eta^n((\na \vel^{n+\frac{1}{2}}+(\na\vel^{n+\frac{1}{2}})^T)\cdot \boldsymbol{n})ds~.
\ee

Multiplying the fourth equation in system \eqref{sys_nd_dis}  by $\frac{\phi_i^{n+1}-\phi_i^{n}}{Re}$ and integration by parts  lead to

\be\label{mu_dis}
\frac{1}{Re}\sum_i\int_\Omega \mu_i^{n+1/2}(\phi_i^{n+1}-\phi_i^{n})d\bx=\frac{\kappa_B}{Re}\sum_i\int_\Omega \frac{1}{2\epsilon}((f_i^{n+1})^2-(f_i^{n})^2)d\bx\nonumber\\
+\frac{\mathcal{M}_s}{Re}\frac{(S(\phi_i^{n+1})-S_{i,0})^{2}-(S(\phi_i^{n})-S_{i,0})^{2}}{2S_{i,0}}+\frac{\alpha}{Re}\sum_i\int_\Omega  (H_i^{n+1}-H_i^n)  d\bx+\frac{\alpha}{Re}\sum_i\int_\Omega (f_w(\phi_i^{n+1})-f_w(\phi_i^n)) d\bx\\
-\frac{\kappa_B}{Re}\sum_i\int_{\partial\Omega_w} \partial_n f_i^{n+\frac{1}{2}}(\phi_i^{n+1}-\phi_i^{n})ds
-\frac{\mathcal{M}_s}{Re}\sum_i\int_{\partial\Omega_w} \frac{S(\phi_i^{n+\frac{1}{2}})-S_{i,0}}{S_{i,0}}\epsilon\partial_n\phi_i^{n+\frac{1}{2}}(\phi_i^{n+1}-\phi_i^{n})ds\nonumber~.
\ee 

Multiplying the third equation in system \eqref{sys_nd_dis}  by $\frac{\mu_i^{n+\frac 1 2}\Delta t}{Re}$  yield
\be\label{phit_dis}
&&\frac{1}{Re}\sum_i\int_\Omega \mu_i^{n+\frac 1 2}(\phi_i^{n+1}-\phi_i^{n})d\bx+\frac{\Delta t}{Re}\sum_i\int_{\Omega}\mu_i^{n+\frac 1 2}(\vel^{n+\frac 1 2}\cdot\na)\phi_i^{n+\frac 1 2}d\bx\nonumber\\
&=&-\frac{\mathcal{M}\Delta t}{Re}\sum_i\int_\Omega (\nabla\mu_i^{n+\frac 1 2})^2d\bx~.
\ee
Multiplying the last equation in system \eqref{sys_nd_dis}  by $\frac{\lambda^{n+\frac 1 2}\Delta t}{Re}$ and integration by parts then sum by $i$ give 

\be\label{lamda_dis}
\frac{\Delta t}{Re}\sum_i\int_\Omega (\lambda_i^{n+\frac{1} 2}\delta_{\epsilon_i} \mathcal{P}_i^n)  :\nabla\vel^{n+\frac 1 2} d\bx-\frac{\Delta t}{Re}\sum_i\int_\Omega \xi \epsilon^2 (\phi_i^n)^2 \left|\nabla \lambda_i^{n+\frac 1 2}  \right|^2 d\bx=0~.
%-\frac{\Delta t}{Re}\sum_i\int_\Omega \xi \epsilon^2 (\phi^n)^2 \left|\nabla \lambda^{n+\frac 1 2}  \right|^2 d\bx
\ee
The discretized energy dissipation law \eqref{energye3_dis} is obtained by combining Eqs.~\eqref{kenetic_dis}-\eqref{lamda_dis} and organizing the terms according to the boundary conditions $L(\phi_i)$ as shown in \eqref{bd_nd_dis}. $\null \hfill \blacksquare$

\subsection{Fully discrete $C^0$ finite element scheme}
The spatial discretization using  $C^0$ finite element  is straight forward. Let $\Omega$ be the domain of interest with a  Lipschitz-continuous boundary $\partial\Omega$. Let ${\mathbf{W}_b}^h \subset \mathbf{W}_b $ be a finite  element space with respect to the triangulation of the domain $\Omega$. The fully discrete scheme of the system is  to find $$\left(\{\Phi_h\}^{n+1}, \{\boldsymbol{U}\}_h^{n+1}, \{F_h\}^{n+1}, \{\Lambda_h\}^{n+1},  \{p_h\}^{n+1},  \{\vel_h\}^{n+1}\right) \in \mathbf{W}_b^h,$$ such that for any $( \psi_{1,h},...,\psi_{N,h},\chi_{1,h},..., \chi_{N,h}, \zeta_{1,h},...,\zeta_{N,h}, \Theta_{1,h},..., \Theta_{N,h}, q_h, \bv_h)\in\mathbf{W}_b^h$, the following scheme holds.
\be
\left\{\begin{array}{ll}
	\label{sys_nd_fully_dis_weak}
	\int_\Omega\left(\frac{\vel_h^{n+1}-\vel_h^{n}}{\Delta t}+(\vel_h^{n+\frac{1}{2}}\cdot\na) \vel_h^{n+\frac{1}{2}}+\frac{1}{Re}\na P_h^{n+\frac{1}{2}}\right)\cdot\bv_h d\bx\\
	=-\int_\Omega\frac{1}{Re}(\eta_h^{n}(\na \vel_h^{n+\frac{1}{2}} +(\na \vel_h^{n+\frac{1}{2}})^T)):\na\bv_h d\bx\\
	~~~~~~~~~~ +\sum_i\int_\Omega\frac{1}{Re}\mu_{i,h}^{n+\frac{1}{2}}\na\phi_h^{n+\frac{1}{2}}\cdot\bv_h d\bx-\sum_i\int_\Omega \frac{1}{Re}\lambda_{i,h}^{n+\frac{1}{2}} \mathcal{P}_{i,h}^{n} \delta_{i,h,\epsilon}:\bv_h d\bx\\
	~~~~~~~~~~ +\int_{\partial\Omega_w}\frac{1}{Re}\bn\cdot(\eta_h^{n}(\na \vel_h^{n+\frac{1}{2}} +(\na \vel_h^{n+\frac{1}{2}})^T)+\sum_i\lambda_{i,h}^{n+\frac{1}{2}} \mathcal{P}_{i,h}^{n} \delta_{i,\epsilon})\cdot\bv_h d\bx~,\\
	\int_\Omega(\na\cdot \vel_h^{n+\frac{1}{2}})q_h d\bx=0~,\\[3mm]
	\int_\Omega(\frac{\phi_{i,h}^{n+1}-\phi_{i,h}^{n}}{\Delta t}+(\vel_h^{n+\frac{1}{2}}\cdot\na)\phi_{i,h}^{n+\frac{1}{2}})\psi_{i,h} d\bx=-\int_\Omega\mathcal{M}\na\mu_{i,h}^{n+\frac{1}{2}}\na\psi_{i,h} d\bx~,\\[3mm]
	\int_\Omega\mu_{i,h}^{n+\frac{1}{2}}\chi_{i,h} d\bx = \int_\Omega\left(    \kappa_B \frac{1}{\epsilon^2}((\phi_{i,h}^{n+1})^2+(\phi_{i,h}^{n})^2+\phi_{i,h}^{n+1}\phi_{i,h}^n-1)f_{i,h}^{n+\frac{1}{2}} \right. \\
	~~~~~~~~~~ \left. +\mathcal{M}_s\frac{(S(\phi_{i,h}^{n+\frac{1}{2}})-S(\phi_{i,h,0}))}{S(\phi_{i,h,0})}(\frac{1}{4\epsilon}((\phi_{i,h}^{n+1})^2+(\phi_{i,h}^{n})^2-2)(\phi_{i,h}^{n+1}+\phi_{i,h}^n))\right)\chi_{i,h} d\bx\\
	~~~~~~~~~~ +\int_\Omega(\kappa_B\na f_{i,h}^{n+\frac{1}{2}}+\mathcal{M}_s\epsilon\frac{(S(\phi_{i,h}^{n+\frac{1}{2}})-S(\phi_{i,h,0}))}{S(\phi_{i,h,0})} \na \phi_{i,h}^{n+\frac{1}{2}})\cdot\na\chi_{i,h} d\bx\\
	~~~~~~~~~~
	+ \int_\Omega \alpha \frac{f_{w}(\phi_{i,h}^{n+1})-f_{w}(\phi_{i,h}^{n})}{\phi_{i,h}^{n+1}-\phi_{i,h}^n} \chi_{i,h} d\bx+ \int_\Omega \alpha \frac{(H_{i,h}^{n+1}-H_{i,h}^n)}{\phi_{i,h}^{n+1}-\phi_{i,h}^n} \chi_{i,h} d\bx\\
	~~~~~~~~~~ -\int_{\partial\Omega_w}(\kappa_B\partial_\bn f_{i,h}^{n+\frac{1}{2}}+\mathcal{M}_s\epsilon\frac{(S(\phi_{i,h}^{n+\frac{1}{2}})-S(\phi_{i,h,0}))}{S(\phi_{i,h})} \partial_\bn \phi_{i,h}^{n+\frac{1}{2}})\chi_{i,h} ds~,  \\ 
	\int_\Omega f_{i,h}^{n+\frac{1}{2}}\zeta_{i,h} = \int_\Omega\epsilon \na\phi_{i,h}^{n+\frac{1}{2}}\cdot\na\zeta_{i,h}+\int_\Omega\frac{1}{\epsilon}((\phi_{i,h}^{n+\frac{1}{2}})^2-1)\phi_{i,h}^{n+\frac{1}{2}}\zeta_{i,h} d\bx\\
	~~~~~~~~~~~-\int_{\partial\Omega_w} \eps\partial_\bn \phi_{i,h}^{n+\frac{1}{2}}\zeta_{i,h} d\bx~,\\
	\int_\Omega	\delta_{i,h,\epsilon} \mathcal{P}_{i,h}^{n} : \na\vel_{i,h}^{n+\frac{1}{2}}\Theta_{i,h} d\bx-\int_\Omega \xi \epsilon^2 ((\phi_{i,h}^{n})^2\nabla \lambda_{i,h}^{n+\frac1 2})\cdot\na\Theta_{i,h}d\bx+\int_{\partial\Omega_w}\xi \epsilon^2 ((\phi_{i,h}^{n})^2\partial_\bn \lambda_{i,h}^{n+\frac1 2})\Theta_{i,h} d\bx=0~.
\end{array}\right.
\ee 

Newton's iteration method is applied to solve the above nonlinear system.    The unique solubility of (\ref{sys_nd_fully_dis_weak}) can be proved by following the idea introduced in \cite{shen2022energy}.

\section{Numerical Results}
\label{sec:results}
In this section,  we first calibrate the model parameters and validate the model by comparing with the experimental data on RBC deformation under different stretching forces.  Then cell-wall interaction simulations are used to study the effects of adhesion force and local insensibility. 
The energy stability of the numerical scheme is also illustrated. Then the calibrated model is used to study the cell-wall interaction and cell-cell aggregation. 

\subsection{BenchMark: vesicle deformation under stretch}

\begin{figure}[!ht]
	\centering
	\includegraphics[width=4.in]{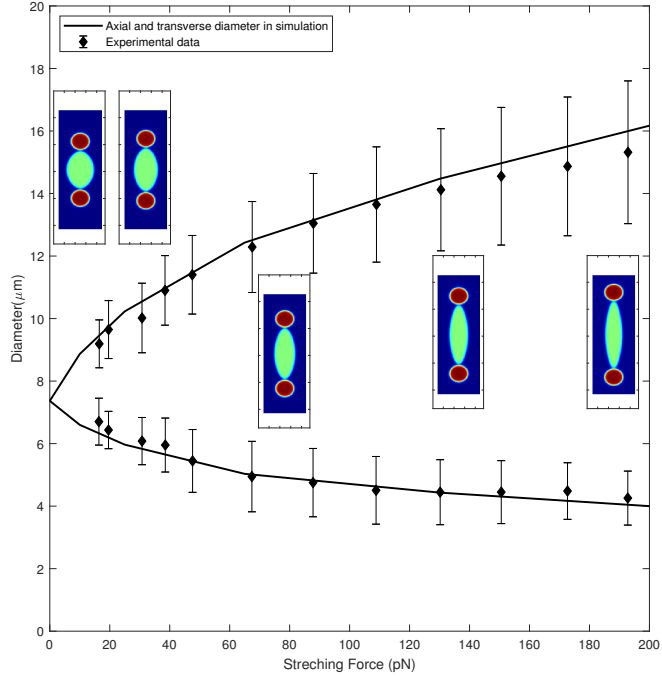}
	\caption{Nonlinear elastic deformation of red blood cell. The curve shows the change of the diameter versus stretching force. The diamond is the experiment data published in \cite{mills2004nonlinear}. In the experiment schematic diagram, the centered phase is set to be the cell and the circulate area on the two side stands for the light tweezers. The force applied on the membrane is kept a constant when the tweezers moving. The equilibrium is reached when the membrane no longer extends under the certain stretch force. }
	\label{fig:diameter}
\end{figure}

% \textcolor{red}{what does "Nonlinear Elastic and Viscoelastic" mean?}
% \textcolor{blue}{Studies on mechanical behaviour of red blood cell shows that the cell's deformation under linearly increased stretching force is highly non-linear \cite{puig2007viscoelasticity,mills2004nonlinear}.This kind of characteristic is named non-linear elasticity or, more precisely , viscoelasticity. } 
Laboratory experiments have tested the non-linear elasticity and deformation of RBC \cite{mills2004nonlinear} in which  optical tweezers are used to provide stretching force to the cells.  We set up a numerical simulation mimicking RBC deformation in the experiment. New phases $\phi_{tw_1}, \phi_{tw_2}$ are introduced to represent optical tweezers. The  energy corresponding to optical tweezers stretching RBC is defined as:
\be
H_{tw_i}= Q_{tw1}  (\phi_1+1)^2(\phi_{tw_i}+1)^2- Q_{tw2} (\phi_1^2-1)^2(\phi_{tw_i}^2-1)^2~, ~~i=1,2~,
\ee
where $Q_{tw}$ is a coefficient that controls the attraction (???) force that the optical tweezers provide. The value of other parameters are listed as follows:
$Re=2\times 10^{-4}, \mathcal{M}=0.25, \kappa_B=2\times 10^{-3}, k=2\times 10^{-12},  l_s=5\times 10^{-3}, \mathcal{M}_s=2.$
% \begin{figure}[!ht]
% 	\centering
% 	\includegraphics[width=4.5in]{figure/benchmark/nonlinear.eps}
% 	\caption{ The centered phase is set to be the cell and the circulate area on the two side stands for the light tweezers. The force applied on the membrane is kept a constant when the tweezers moving. The equilibrium is reached when the membrane no longer extends under the certain stretch force.}
% 	\label{fig:bench mark}
% \end{figure}
The force  applied on the cell is calculated by the following equation:
\be
\bF=\int_{\Omega}\frac {\partial H_{tw_i}} {\partial \phi}\na \phi d\bx~.
\ee
The curves of axial and transverse diameter versus stretching forces are shown in Figure \ref{fig:diameter} together with the experimental data (diamonds with bars) from \cite{mills2004nonlinear}.
The results show that our model   fits the experimental data very well. 
%However, the error still exist since the simulation is in 2-D domain and the initial shape of the red blood cell a circle instead of a biconcave curve.

\subsection{Cell-wall attraction}
\begin{figure}[!ht]
	\centering
	\includegraphics[width=4.0in]{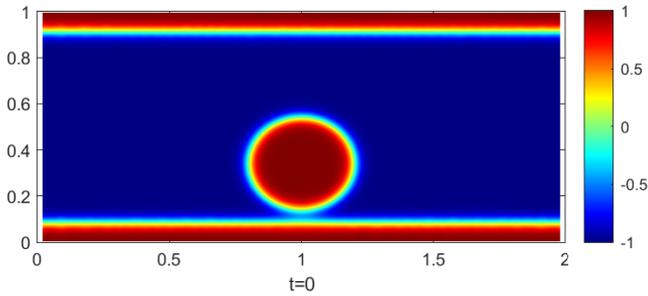}
	\caption{The initial state of the case.}
	\label{fig:init}
\end{figure}

Cell-wall interactions under blood flow conditions  play important roles during blood clotting  \cite{wu2014three} and cancer cell invasion \cite{gerisch2008mathematical}. 
This example is used to  investigate the effect of adhesion force on  the vesicle-wall interaction. We first consider the case without external flow. 
As shown in Figure \ref{fig:init}, a vesicle is initially placed at a location with a point-wise contact with the wall phase. The parameter values of this simulation are listed as follows: 
$Re=2\times 10^{-4},\mathcal{M}=5\times10^{-4},\kappa_B=2\times 10^{-2},\eps=2\times10^{-3},\mathcal{M}_s=10^{2},k=4\times 10^{-11}, l_s=5\times10^{-6}, \alpha=1000, q_{w_1}=2, q_{w_2}=1$.

\begin{figure}[!ht]
	\centering
	\includegraphics[width=5.5in]{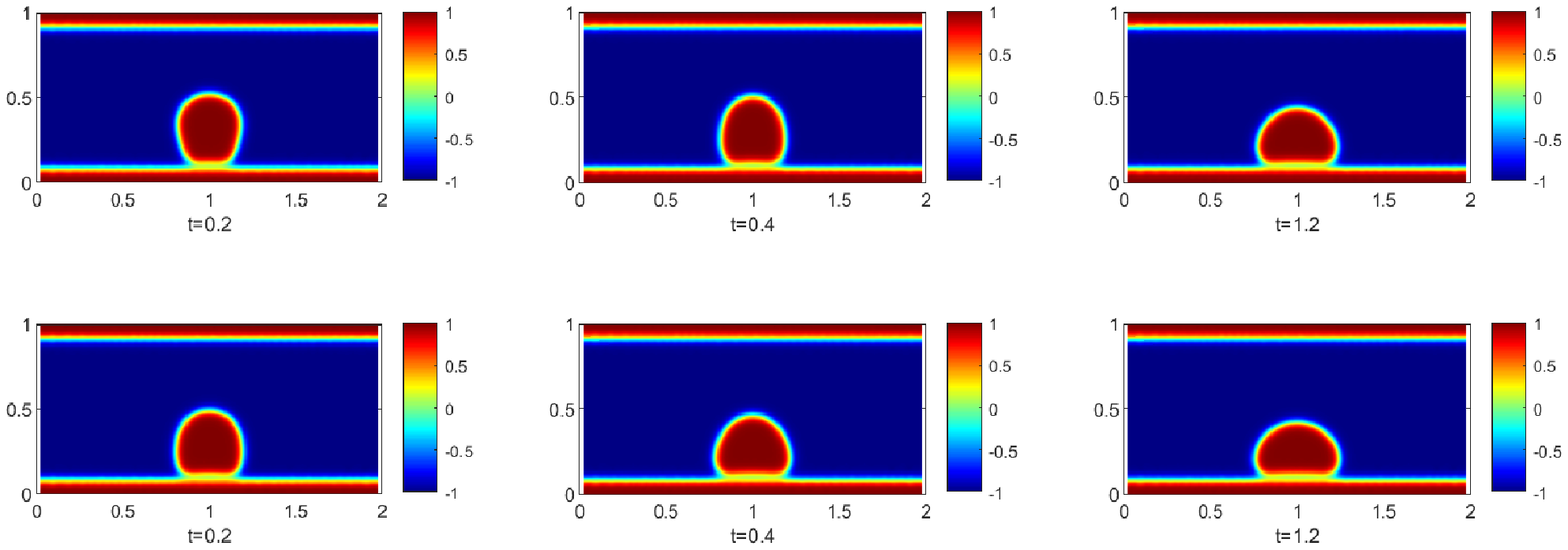}
	\caption{The top three pictures shows  the deformation of the vesicle with local inextensibility. The bottom three pictures shows the deformation of the vesicle without local inextensibility.}
	\label{fig:wall attraction}
\end{figure}

\begin{figure}[!ht]
	\centering
	\includegraphics[width=5.5in]{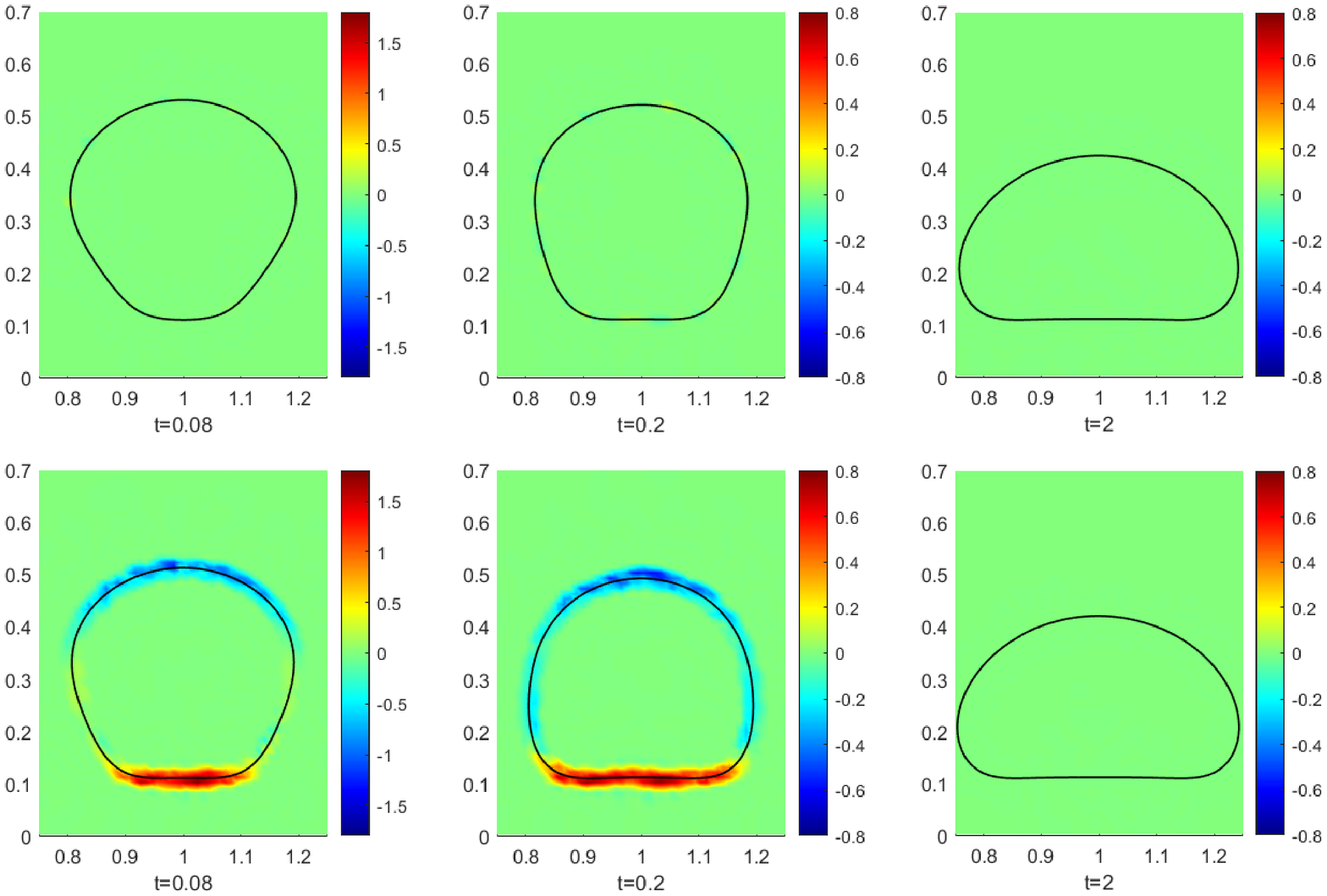}
	\caption{The top two pictures shows  the value of surface divergence $\mathcal{P}_i: \na\vel$ with local inextensibility. The bottom pictures shows the surface divergence of the vesicle without local inextensibility.}
	\label{fig:surfacediv}
\end{figure}

\begin{figure}[!ht]
	\centering
	\includegraphics[width=3.5in]{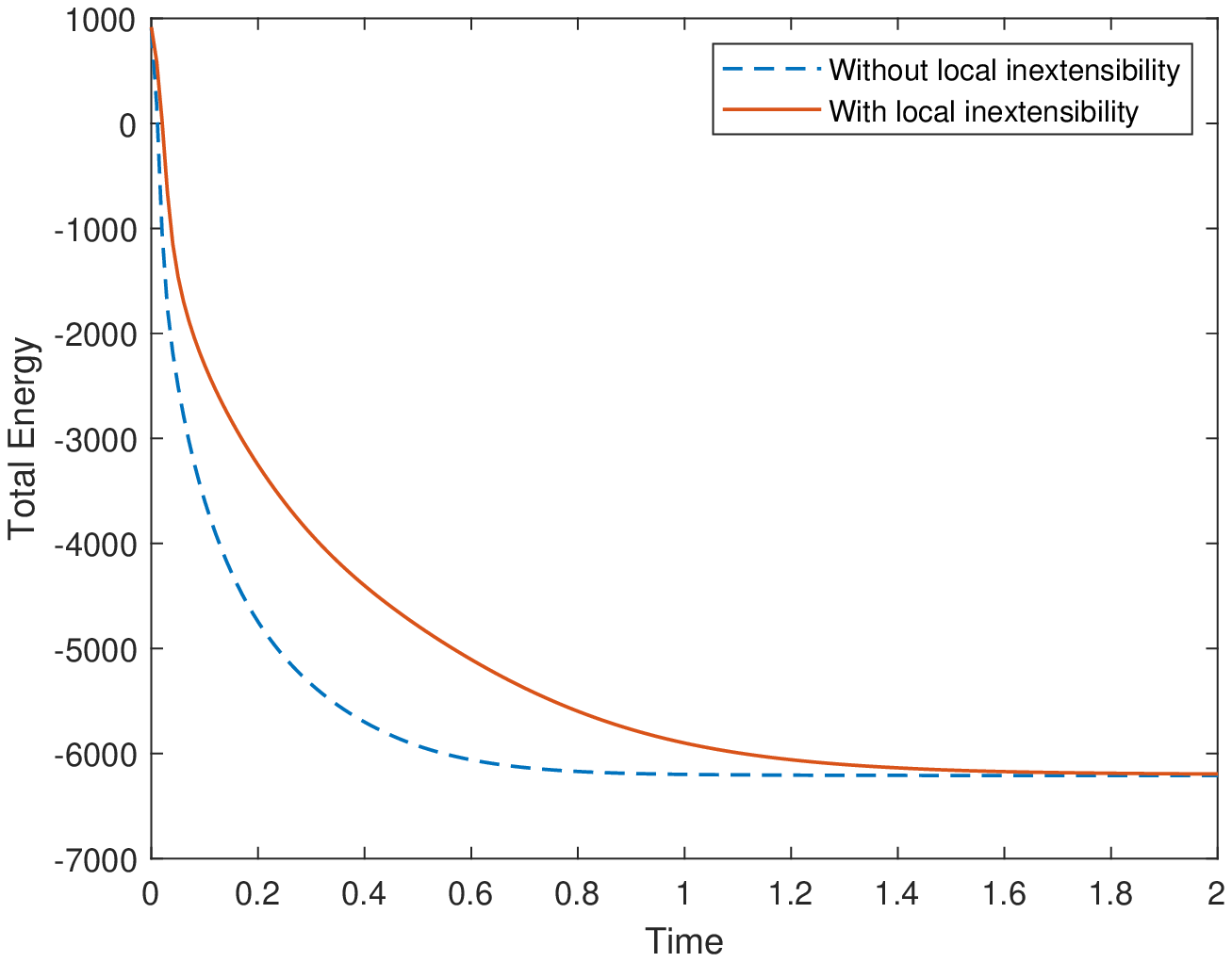}
	\caption{Total energy of the two system versus time.}
	\label{fig:energy}
\end{figure}

Figure \ref{fig:wall attraction} shows the scenarios of the the vesicle with and without local inextensibility at different time.   As shown in Fig. \ref{fig:surfacediv}, the local inextensiblity constrains slows down the deformation of the vesicle.  For vesicle without local constrain, it allowed the membrane attached to the wall to be stretched (red color in Fig.\ref{fig:surfacediv}) to achieve the equilibrium faster. This also could be seen in the Fig. \ref{fig:energy}. It confirms that in both cases, the energy monotonously decays to the same value while the vesicle without local inextensibility evolves faster.

\begin{figure}[!ht]
	\centering
	\includegraphics[width=5.5in]{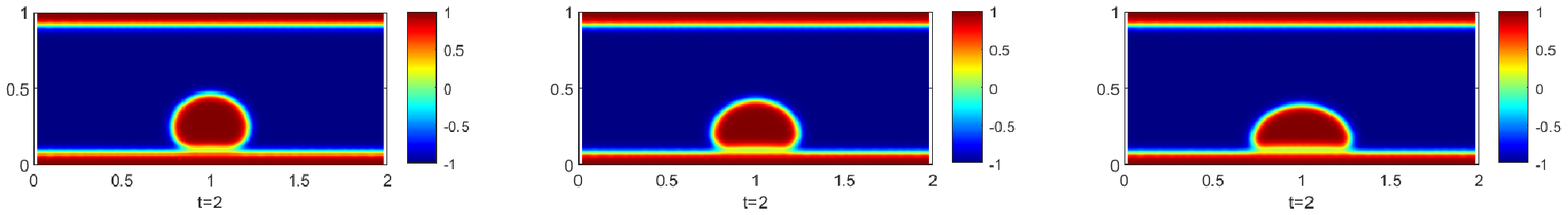}
	\caption{Equilibrium of the cells at different strength of the adhesion force. (Left: $\alpha=400$; Mid: $\alpha=1000$; Right: $\alpha=2500$.)}
	\label{fig:diffstrength}
\end{figure}

Then the effect of the strength of the adhesion force  on the equilibrium profiles of the vesicles are illustrated in Figure \ref{fig:diffstrength}.  As expected, when the strength increases from 400 to 2500,  the contact surface area increases. 

\begin{figure}[!ht]
	\centering
	\includegraphics[width=5.5in]{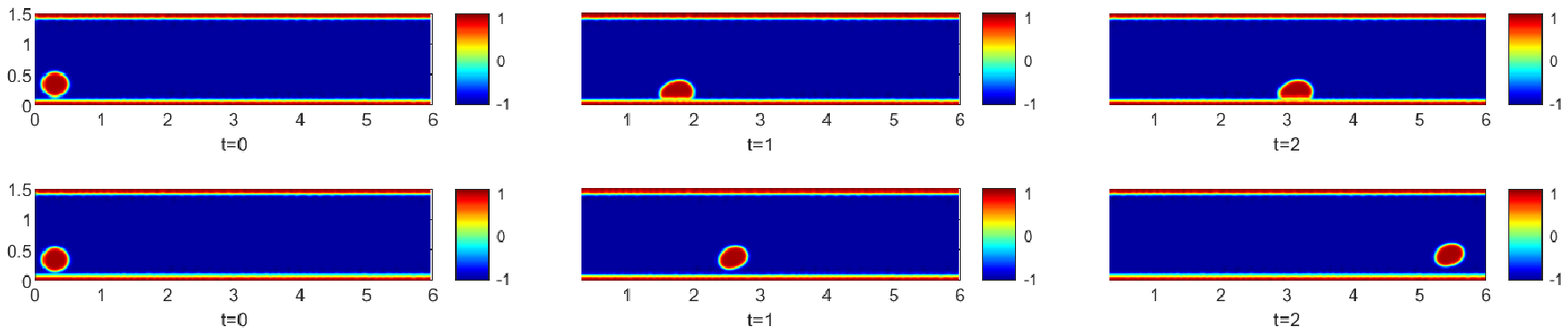}
	\caption{The top three picture shows the motion of the cell in strong adhesion case. The bottom three pictures shows the cell motion in weak adhesion case.}
	\label{fig:shear}
\end{figure}

During the blood clotting, the platelets interact  with several coagulation factors and get activated. The full activated plated has strong adhesion force with vessel wall in order to stop bleeding more efficiently. 
In Figure \ref{fig:shear},  we  simulate two platelets under the shear blood flow, one fully activated with strong adhesion  and one partial activated with weak adhesion with vessel wall.   It shows that when the fully activated one is captured by the wall; while the partial activated one is washed away by the shear flow. 

% The corresponded surface area and volume versus time is also shown in Figure \ref{fig:surf_vol}

% \begin{figure}[!ht]
% 	\centering
% 	\includegraphics[width=3.5in]{figure/wallattraction/surfacevolume.eps}
% 	\caption{The curves show how the surface and volume change with time.}
% 	\label{fig:surf_vol}
% \end{figure}

\subsection{Cell aggregation and offset}
\begin{figure}[!ht]
	\centering
	\includegraphics[width=5.5in]{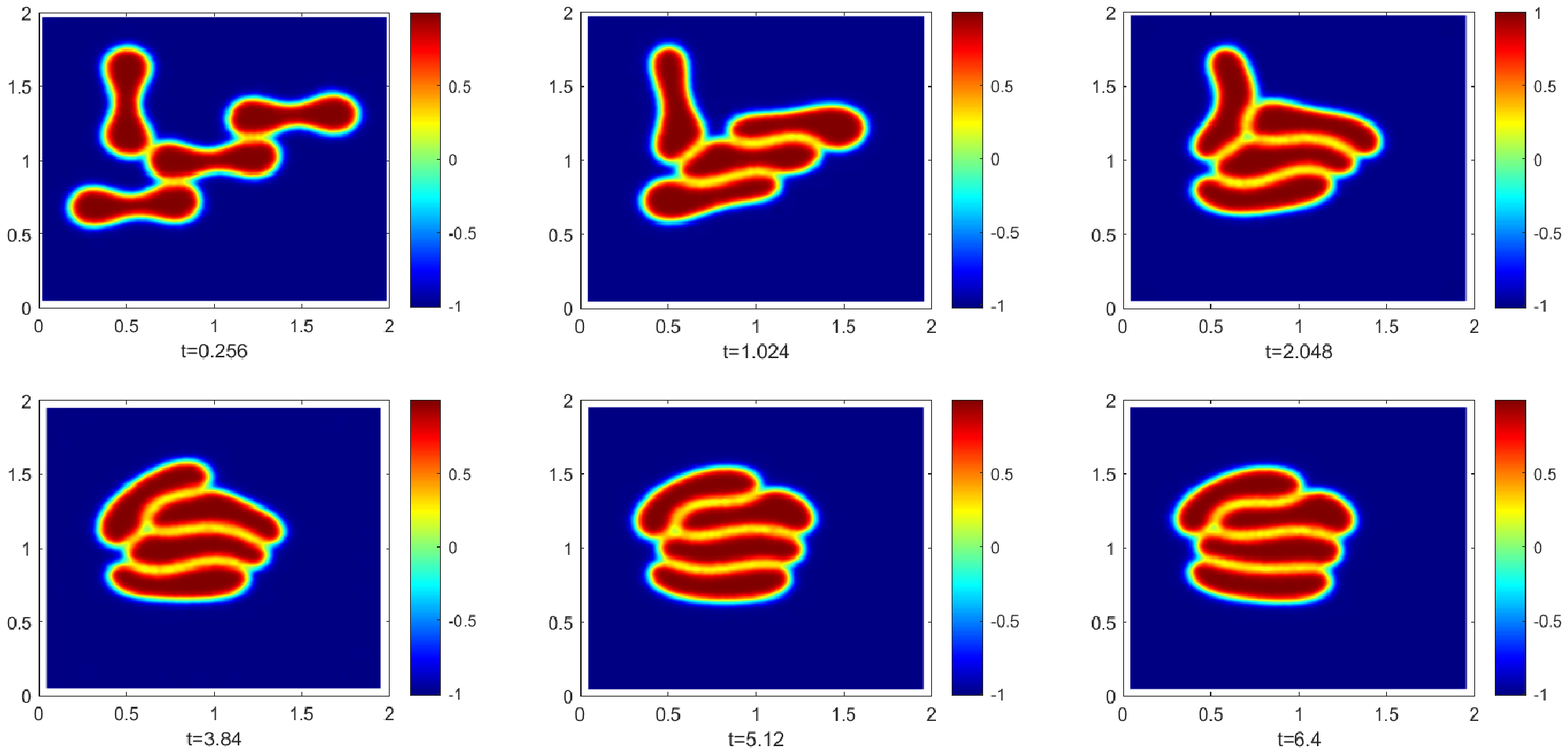}
	\caption{Aggregating of four red blood cells.}
	\label{fig:aggregating}
\end{figure}

Aggregation of red blood cell is a phenomenon observed in experiment \cite{Thomas2008aggregate}. In this part, we set up a simulation with four red blood cells contacted with each other at a small point. The parameters are shown below.
$Re=2\times 10^{-5},\mathcal{M}=5\times 10^{-4}, \kappa_B=4\times 10^{-2}, k=4\times 10^{-12}, l_s=5\times 10^{-3}, \mathcal{M}_s=10^3, \alpha=3\times10^3,q_1=1, q_2=0.5$
The evolution of the system with time is shown in \ref{fig:aggregating}.

From the result we can see that they are creeping together with time under the attractive force and form rouleaux which is consistent with the experimental result shown in \cite{Thomas2008aggregate} and \cite{ziherl2007aggregates}.

Also, the deformation of the cells at equilibrium is related to the value of the attractive force \cite{ziherl2007aggregates,ziherl2007flat,flormann2017buckling}. Fig. \ref{fig:offset} shows the status of moderately and strongly aggregate. The parameters are shown below.
$Re=2\times 10^{-4}, \mathcal{M}=5\times 10^{-4}, \kappa_B=2\times 10^{-2}, k=2\times 10^{-11}, l_s=5\times 10^{-3},  \mathcal{M}_s=10.$

\begin{figure}[!ht]
	\centering
	\includegraphics[width=4.5in]{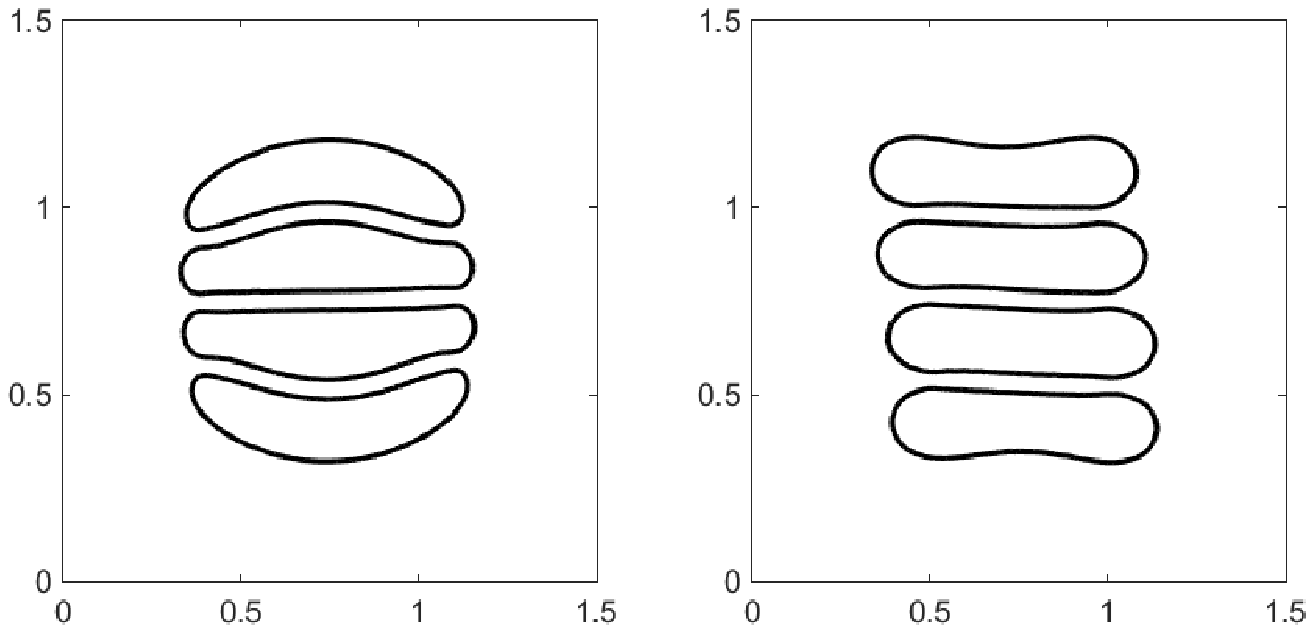}
	\caption{Left: Strong aggregate $\alpha=15\times10^3, q_1=1, q_2=1$. Right: moderate aggregate $\alpha=100, q_1=1, q_2=1$. }
	\label{fig:offset}
\end{figure}
\begin{figure}[!ht]
	\centering
	\includegraphics[width=6.5in]{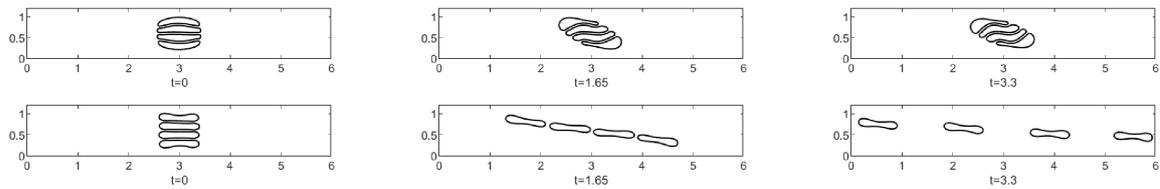}
	\caption{The top three figures shows the motion of the cells with strong aggregation. The bottom three figures shows the motion under moderate aggregation. }
	\label{fig:couette}
\end{figure}

From the figure we can see that under strong aggregation force, obvious terminal
hemispherical caps is shown and the offset between each adjacent cell is smaller as well. This result is consistent with the experimental phenomenon observed well \cite{Thomas2008aggregate}.
In the following test we set the cells on a Couette flow with shear rate equal to be  $ 20s^{-1}$ with dimension as in \cite{zhang2008red}. The motion of the cells are shown in figure \ref{fig:couette}.
Under shear flow, the rouleaux with strong adhesion still  aggregate together; while the one with weak adhesion is broken up by the shear flow, which is consistent with the  results reported in \cite{zhang2008red}.

%From the figure we can see that under strong aggregation force the thickness of the cells are smaller and the offset between each adjacent cell is smaller as well. This result fits the experimental phenomenon observed well. Also, by varying the coefficients that governs the attractive and repulsive force, the overlapping extent of different phases at equilibrium, i.e. the distance between cells can be adjusted manually.

\subsection{Red blood cell motion at vessel bifurcation}
In this part, we simulate the motion of RBCs in branched vessels. 4 cells are initially placed in the Y shaped vessel. Width of the main channel of the Y shaped vessel and the bottom branch of the vessel is set to be $1\times 10^{-7}$ meter. The top branch of the vessel is $0.7\times 10^{-7}$ meter, which is close to the size of a red blood cell. A pressure drop boundary condition is used to  introduce a shear flow in the vessel with the velocity around $5\times 10^{-4}m/s$, which is close to the blood flow in capillaries. Other model parameters in this simulation are as follows:
$Re=2\times 10^{-4}, \mathcal{M}=5\times 10^{-4}, \kappa_B=4\times 10^{-2}, k=4\times 10^{-11}, l_s=2=5\times 10^{-6}, \mathcal{M}_s=20$.

\begin{figure}[!ht]
	\centering
	\includegraphics[width=4.5in]{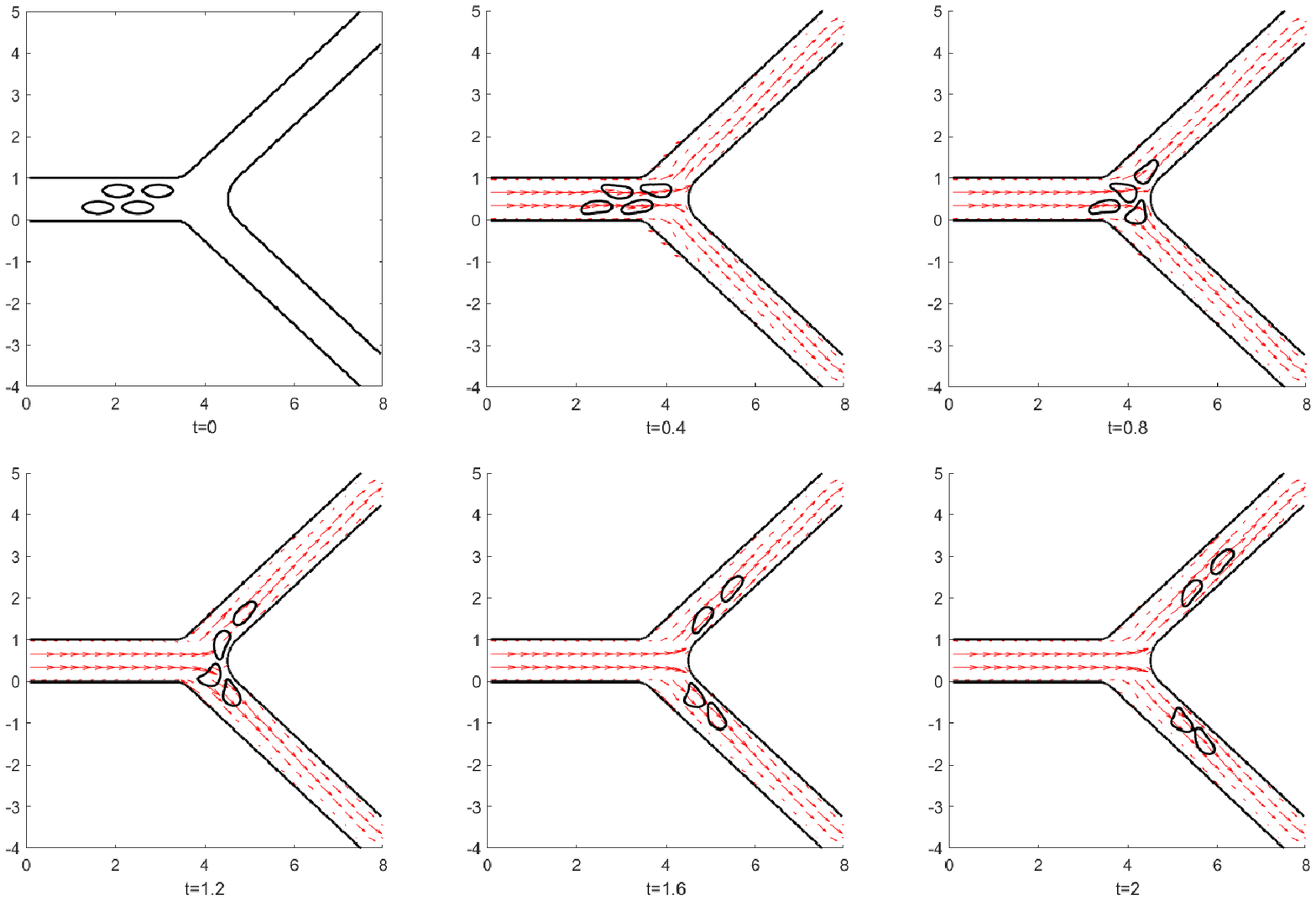}
	\caption{Cells are set in a cluster initially under a moderate aggregating force with $\alpha=25, q_1=1, q_2=1$. The main channel width is 1 and the width of the two branches  is 0.7. RBCs divide equally at the vessel bifurcation. The velocity field is indicated by the vector field.}
	\label{fig:Yinitial}
\end{figure}

\begin{figure}[!ht]
	\centering
	\includegraphics[width=4.5in]{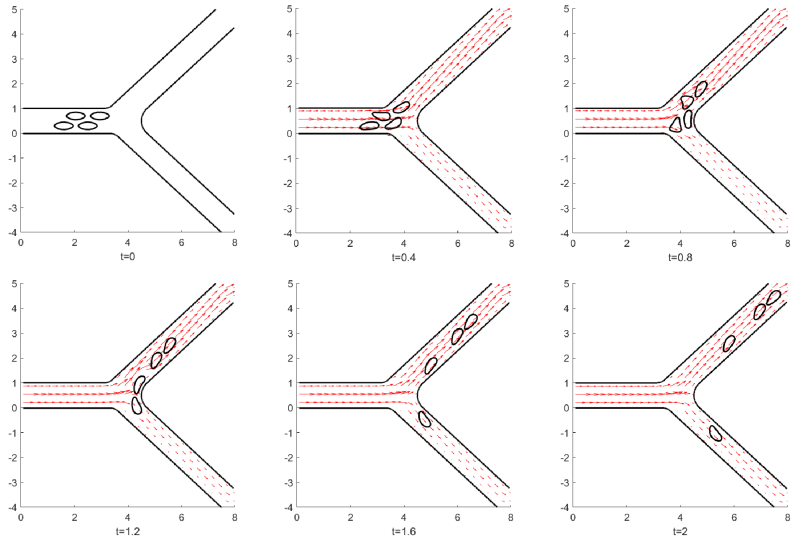}
	\caption{Cells are set in a cluster initially under a moderate aggregating force with $\alpha=25, q_1=1, q_2=1$. The top branch is set to be 1 with the other stays 0.7. One of the cells is going into the branched vessel. The velocity field is shown as well.}
	\label{fig:Yweak}
\end{figure}

\begin{figure}[!ht]
	\centering
	\includegraphics[width=4.5in]{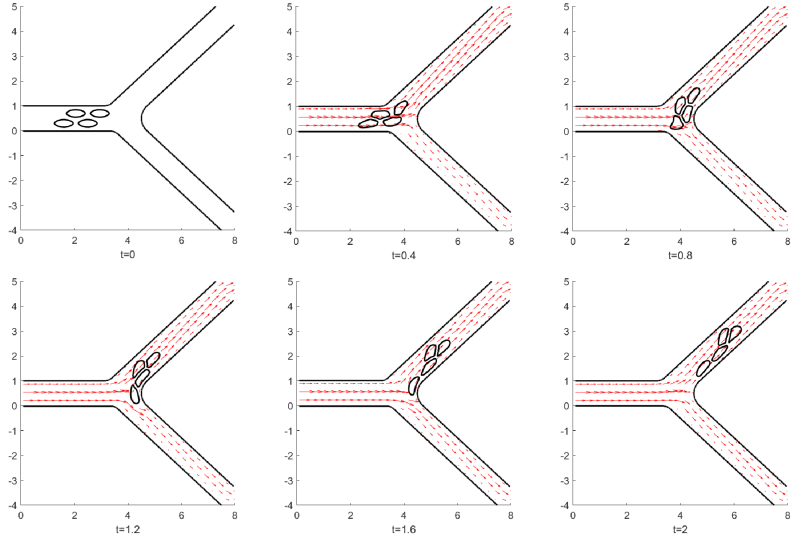}
	\caption{Cells are set in a cluster initially under a strong aggregating force with $\alpha=1.5\times 10^3, 1_1=1, Q_2=1$. The set up is the same as Figure \ref{fig:Yweak}. None of the cells is going into the branched vessel. The velocity field is shown as well.}
	\label{fig:Ystrong}
\end{figure}

Motion of the RBCs and the velocity field of the flow are shown in Figure \ref{fig:Yweak}, \ref{fig:Ystrong}, and \ref{fig:Yinitial}, respectively.
We firstly simulate motion of the RBC group with a moderate aggregation force when they pass a Y-shaped channel with the same width of the both branches as a baseline. The cells divide equally at the vessel bifurcation. 
Then one of the channel is widen and a non-equal deviation is observed (3 move upwards,1 moves downward) due to the  lower resistance in the upper branch. Then,  under the same geometric setup, a strong aggregation force is applied to the cells. In this case,    all of the cells goes into the wider channel.
The simulation result explains the experimental result in \cite{Thomas2008aggregate}, which a great number of red blood cells are observed to be absent in the branched vessel under strong aggregate case compared with moderate aggregate.

%\begin{figure}[!ht]
%	\centering
%	\includegraphics[width=4.5in]{figure/bifurcation/weak2.eps}
%	\caption{Cells are set side by side under a moderate aggregating force with $Q_1=1000, Q_2=250$. Two of the cells are observed to go into the branched vessel.}
%	\label{fig:weak2}
%\end{figure}

%\begin{figure}[!ht]
%	\centering
%	\includegraphics[width=4.5in]{figure/bifurcation/strong2.eps}
%	\caption{Cells are set side by side under a moderate aggregating force with $Q_1=1000, Q_2=600$. None of the cells is observed to go into the branched vessel.}
%	\label{fig:strong2}
%\end{figure}

\section{Conclusions}
\label{sec:conclusion}
In this paper, a thermodynamically consistent diffuse interface (phase-field) model for describing cell-wall, cell-cell interaction and  aggregation is derived based on the energy variation method.  The interactions of cellular and wall structures are modeled  by introducing a new interaction energy with the help of the phase-field function.  

Then an efficient scheme using $C^0$ finite element spatial discretization and the mid-point temporal discretization is proposed to solve the obtained model equations. Thanks to the mid-point temporal discretization, the obtained numerical scheme is unconditionally energy stable.  The model and parameters are calibrated with experimental data on cell deformation under different stretch force. Then effects of  the  adhesion strength  on cell-wall and cell-cell interaction are studied. In the end, the calibrated model is used to model red blood cell motion near  the vessel  bifurcation. The results are consistent with the experimental observation.

\newpage
\bibliographystyle{elsarticle-num} 
\bibliography{Mybib}

\end{document}